\documentclass[12pt]{article}
\usepackage{amsmath}
\usepackage{graphicx}
\usepackage{enumerate}
\usepackage{mathrsfs,accents,amssymb,multirow,longtable,rotating,array,booktabs,float,bm}
\usepackage{natbib}
\usepackage{hyperref}
\usepackage{pdflscape}
\usepackage{url} 
\usepackage{color}
\usepackage{wasysym}
\usepackage{arydshln}
\newcommand{\blind}{1}

\newcommand{\ba}{\begin{equation}}
\newcommand{\ea}{\end{equation}}

\newtheorem{theorem}{Theorem}
\newtheorem{proposition}{Proposition}

\newtheorem{remark}{Remark}

\allowdisplaybreaks
\begin{document}

\def\spacingset#1{\renewcommand{\baselinestretch}%
{#1}\small\normalsize} \spacingset{1}


\if1\blind
{\title{\bf Matrix-variate integer-valued autoregressive processes}
  \author{Nuo Xu\\
Jilin University,  China,\\
Changchun University of Technology, China\\
    Kai Yang\thanks{
    Corresponding author. E-mail: yangkai@ccut.edu.cn. Yang's research was supported by National Natural Science Foundation of China (No. 12471249) and Natural Science Foundation of Jilin Province (No. 20220101038JC). Zhu’s research was supported by National Natural Science Foundation of China (No.
12271206) and Natural Science Foundation of Jilin Province (No. 20250102003JC).} \\
Changchun University of Technology, China\\
    and \\
    Fukang Zhu \\
  Jilin University,  China}
   \maketitle
} \fi

\if0\blind
{
  \bigskip
  \bigskip
  \bigskip
  \begin{center}
    {\LARGE\bf Matrix-variate integer-valued autoregressive processes}
\end{center}
  \medskip
} \fi

\bigskip
\begin{abstract}
In the fields of sociology and economics, the modeling of matrix-variate integer-valued time series is urgent. However, no prior studies have addressed the modeling of such data. To address this topic, this paper proposes a novel  matrix-variate integer-valued autoregressive model. 
The key techniques lie in defining two left- and right-matricial thinning operators. 
The probabilistic and statistical properties of the proposed model are investigated. 
Furthermore, two estimation methods are developed: projection estimation and iterative least squares estimation.  The corresponding asymptotic properties of these estimators are established. Additionally, the order-determination problem is addressed. In the simulation studies, the estimation results are given and the theoretical properties are verified. Finally, it is shown that the matrix-variate integer-valued autoregressive model is superior to the continuous matrix-variate autoregressive and multivariate integer-valued autoregressive models for  matrix-variate integer-valued time series data.
\end{abstract}

\noindent%

\vfill

\spacingset{1.9}
\section{Introduction}

In recent years, with the rapid development of big data, artificial intelligence and other fields,
data collection and storage technology has been continuously improved,
and the structure of observed data has become increasingly complex.
In various fields such as economy, finance, medicine, and insurance, not only vector or high-dimensional data (\cite{Lam2012,Wang2022,Zhang2024}) 
 formed by observing different indicators of the same phenomenon,
but also  panel data (\cite{Jiang2023}) obtained by observing the same indicator in different regions have emerged. 
These two types of data are usually highly correlated. Analyzing one of them alone will cause serious information loss
and cannot meet the needs of modern data analysis.
Therefore, people try to observe these two classifications together, thus forming a more complex matrix-variate data.

Compared with  point-valued or vector-valued data, 
matrix-valued data usually have a more complex structure and contain richer information. In order to better understand matrix-variate data, we take   two examples in Section A of Supplement Material.
Analyzing matrix-variate data is an emerging challenge in data analysis. 
First, the dependency structure between rows and columns of matrix-variate data cannot be replaced by traditional vector data; Second, the increase in data dimensions may led to the complexity of the model structure, which poses great challenges to parameter estimation. 
Early studies on matrix-variate data focused on the distribution types of matrix variables (\cite{Rukhin2003,Chikuse1981,Dawid1981}.
At the beginning of this century, scholars began to consider the matrix-variate modeling.  
\cite{Wang2009, Leng2012} considered statistical analysis for different matrix normal graphical models. 
\cite{Zhou2014}  proposed a class of regular matrix regression methods based on spectral regularization. 
\cite{Ding2018, Kong2020}   studied the regression analysis of matrix-variate response under different situations. 
Some recent achievements about this topic refer to 
\cite{Chang2023, Han2024}, among others.

The matrix-variate data are usually observed by time, thus forming matrix-variate time series (MTS). 
If we convert MTS into vector time series for modeling, 
it will destroy the matrix structure of the data and sever the cross relationship between rows and columns, 
thereby cause serious information loss.
Currently, the methods for MTS modeling can be divided into two categories: 
The first one is matrix factor model (\cite{Wang2019,Chen2020,Yu2022,Chen2023b,Yuan2023}). 
This kind of model can maintain the structure of matrix data and effectively reduce the number of parameters that need to be estimated. Therefore, it can well deal with the modeling problem of high-dimensional MTS data. 
The other type is the matrix autoregressive (MAR) model (\cite{Chen2021}), which can well describe the dependence of rows and columns of MTS data through left and right multiplication coefficient matrices. 
Following the MAR model, \cite{Hsu2021}  and \cite{Zhang2023} successively proposed the spatio-temporal MAR model and the additive MAR model, and considered the application of the corresponding models in different fields. 
\cite{Tsay2024}  reviewed the recent achievements on MTS,  
and further proposed the matrix-variate autoregressive moving average model. 
Recently, \cite{Yu2024} proposed a novel matrix generalized autoregressive conditional heteroscedasticity  for matrix-variate financial series.  
These results have made great contributions for MTS.

As is mentioned previously, 
in addition to the continuous MTS data mentioned above, there are also a large number of 
MTS taking values on non-negative integers in our daily life. 
There is strong evidence to prove that the word frequency statistics of sensitive vocabulary can effectively monitor the development trend of online public opinion.
Also, the criminal activities data in different regions play an important role in promoting the correlation analysis of criminal activities and the detection of cases. 
However, a continuous MTS model may not suitable for fitting such data. 
Not only is it  difficult to provide integer predictions, but the model is also not very interpretable.
Therefore, it is necessary to propose a reasonable statistical model that specifically for  matrix-variate integer-valued time series (MITS) data.

A popular method for modeling integer-valued time series data is based on the tinning operators (\cite{Steutel1979}). 
The thinning 
operations are applied to random counts, and always lead to integer-valued results. 
This makes the thinning models very interpretable when modeling dependent count data. 
\cite{Alosh1987}  pioneered the first-order integer-valued autoregressive (INAR) model based on the binomial thinning operator, marking the birth of integer-valued time series analysis. 
In recent years, 
with the continuous deepening of theoretical research and the development and change of application scenarios, 
the study of integer-valued time series has gradually developed into fields such as 
heteroscedasticity (\cite{Fokianos2009}), 
high-dimensional (\cite{Yang2023,Xu2025}), 
and 
nonlinear scenarios (\cite{Wang2014,Lih2024}). 

However, 
the research on MITS modeling is still an open problem. 
To model MITS, 
there are at least two critical issues  that need to be addressed. 
The first problem is the modeling mechanism of the matrix-variate count data.  
Traditional thinning operators are not applicable. 
Therefore, some kind of new thinning operations specifically for MITS  need to be defined first. 
The second problem is how to control the number of model parameters, 
because too many parameters will lead to a decrease in estimation accuracy or even estimation failure.  
Inspired by matrix algebra theory, 
we define the left- and right-matricial thinning operators. 
Under this framework, 
the original matrix structure of the data is fully taken into account, thereby, 
an integer-valued version of the MAR model can be nicely constructed with naturally interpretability.  
In addition, 
it reduces the number of model parameters from   $m^2n^2+mn$ to  $m^2+n^2+mn$ relative to the traditional multivariate INAR models of order one, thus making the estimation method more effective. 
Furthermore, this modeling mechanism can be easily extended to other  matrix-variate integer-valued data analysis scenarios and integrated with various thinning operations, 
which is a considerable improvement on data analysis for time series of counts. 

The rest of the paper is organized as follows. 
In Section \ref{sec:2}, we  introduce the  $p$th order matrix-variate integer-valued  autoregressive (MAT-INAR($p$)) model based on two  new matricial thinning operators, along with some of its properties. 
In Section \ref{sec:3}, projection estimation and the  asymptotic theory of the estimators are given. 
The iterated conditional least squares estimation and the related asymptotic theory are presented in Section \ref{sec:4}.
Section \ref{sec:5} addresses the order determination problem. 
Some numerical studies are carried out in Section \ref{sec:6}, 
while  a practical application example is given in Section \ref{sec:7}. 
All proofs, additional
discussions, tables, figures are collected in Supplement Material.

\section{Matrix-variate integer-valued autoregressive  process}\label{sec:2}
\subsection{Matricial thinning operators}\label{sec:2.1}

In order to model MITS, it is necessary to introduce the definition and properties of matricial thinning operators. Specifically, we introduce the left-matricial  thinning operator ``$\bm{A}{\circ}_L$'' and
the right-matricial thinning operator  ``${\circ}_R\bm{B}$'' in this study.
Let $\bm{A}=(a_{i,j})_{m\times m}$ be a square matrix   with elements $a_{i,j}\in [0, 1]$ and  $\bm{B}=(b_{i,j})_{n\times n}$ be a square matrix  with elements $b_{i,j}\in [0, 1]$.
Denote by $\bm{Y}=(y_{i,j})_{m\times n}$  an  integer-valued matrix variate with $y_{i,j} \in \mathbb{N}_0$.
Then, we define the $(i,j)$ elements of $\bm{A}{\circ}_L\bm{Y}$ and $\bm{Y}{\circ}_R\bm{B}$ by
\begin{equation}\label{circ}
(\bm{A}{\circ}_L\bm{Y})_{i,j}=\sum_{k=1}^m a_{i,k}\circ y_{k,j}, ~(\bm{Y}{\circ}_R\bm{B})_{i,j}=\sum_{k=1}^n b_{k,i}\circ y_{j,k},
\end{equation}
where
\begin{equation}\label{circelement}
a_{i,k}\circ y_{k,j}=
\begin{cases}
\sum_{s=1}^{y_{k,j}} B_{a_{i,k}}^{(s)},\mbox{~if~}y_{k,j}>0,\\
0,\hspace{0.6in}\mbox{~if~}y_{k,j}=0,
\end{cases}
~b_{k,i}\circ y_{j,k}=
\begin{cases}
\sum_{s=1}^{y_{j,k}} B_{b_{k,i}}^{(s)},\mbox{~if~}y_{j,k}>0,\\
0,\hspace{0.6in}\mbox{~if~}y_{j,k}=0,
\end{cases}
\end{equation}
``$\circ$'' is the binomial thinning operator introduced by \cite{Steutel1979}, $\{B_{a_{i,k}}^{(s)}\}$ and $\{B_{b_{k,i}}^{(s)}\}$ are sequences of independent and identically distributed (i.i.d.) Bernoulli random variables such that $P(B_{a_{i,k}}^{(s)}=1)=a_{i,k}=1-P(B_{a_{i,k}}^{(s)}=0)$,
$a_{i,k}\in[0, 1]$,
$P(B_{b_{k,i}}^{(s)}=1)=b_{k,i}=1-P(B_{b_{k,i}}^{(s)}=0)$,
$b_{k,i}\in[0, 1]$.
As can be seen from (\ref{circ}) and (\ref{circelement}) that the left- and right-matricial thinning operators
maintain a similar rule of matrix left and right multiplication operations
while replacing scaler multiplication with the binomial thinning operator during element-wise calculation.
Related properties about the  new matricial thinning operators are shown in Section B.1 of  Supplement Material.

\subsection{Definition of MAT-INAR($p$) process}
We introduce the MAT-INAR($p$) process in the following recursive equation:
 \begin{equation}\label{Mymodelp}
\bm{Y}_t = \bm{A}_1{ \circ}_{L} \bm{Y}_{t-1}{ \circ}_{R}\bm{B}_1^{\top}+\cdots+\bm{A}_p{ \circ}_{L} \bm{Y}_{t-p}{ \circ}_{R}\bm{B}_p^{\top}+ \bm{\mathcal{E}}_{t} ,~t \in \mathbb{Z},
\end{equation}
where
\begin{enumerate}[(i)]
 \setlength{\itemsep}{0pt}
 \setlength{\parskip}{0pt}
 \setlength{\parsep}{0pt}
\item $\bm{Y}_t$ is the $m\times n$ matrix-variate integer-valued observation at time $t,~t\in\{1,...,T\}$;
\item $\bm{A}_l=(a_{i,j}^{(l)})_{m \times m}$ is the $m\times m$ autoregressive coefficient matrix with elements $a_{i,j}^{(l)} \in [0,1]$, while
  $\bm{B}_l=(b_{s,k}^{(l)})_{n\times n}$ is the $n\times n$ autoregressive coefficient matrix with elements $(b_{s,k}^{(l)})_{n\times n}\in [0,1]$,
    $i,j\in\{1,...,m\}$,     $s,k\in\{1,...,n\}$, $l\in\{1,...,p\}$;
\item ``$\bm{A}_l{\circ}_L$'' and  ``${\circ}_R\bm{B}_l$''  ($l\in\{1,...,p\}$) are the matricial thinning operators
      defined in (\ref{circ}) and (\ref{circelement}), respectively;
\item $\{\bm{\mathcal{E}}_t\}$ is a series of  i.i.d. $\mathbb{N}^{m\times n}$-valued random error matrix
      following some discrete matrix-valued distribution with mean matrix $\bm{\Lambda}$,
      $E(\bm{\mathcal{E}}_t\bm{\mathcal{E}}_t^{\top}) < \infty$;
\item For fixed $t$, $\bm{\mathcal{E}}_t$ is assumed to be independent of counting series in $\bm{A}_l{\circ}_L\bm{Y}_{t-l}{\circ}_R\bm{B}_l^{\top}$  and $\bm{Y}_{t-s}$ for all $s\geq 1$. And  $\bm{A}_l{\circ}_L\bm{Y}_{t-l}{\circ}_R\bm{B}_l^{\top}$    are also independent of $\bm{Y}_{t-s}$ for all $s\geq 1$.
\end{enumerate}

\begin{remark}
It follows by {\rm (\text{\bf P}.B2)} in Section B.1 of Supplement Material that   the $\bm{A}_l{ \circ}_{L} \bm{Y}_{t-l}{ \circ}_{R}\bm{B}_l^{\top}$ remains unchanged if  $\bm{A}_l$ and $\bm{B}_l$ are divided or multiplied, respectively, by a same nonzero constant.
Therefore, in order to ensure the identifiability of  Model {\rm (\ref{Mymodelp})},
we use the convention that $\bm{A}_l~ {(l\in \{1,...,p\})}$  is normalized so that its Frobenius norm is one, ${\Vert\bm{A}_l\Vert}_F=1$.
\end{remark}

Let ${\rm vec}(\cdot)$ be the vectorization of a matrix by stacking its columns.
Thus, the MAT-INAR($p$) model defined in (\ref{Mymodelp}) can be rewritten as:
\begin{equation}\label{Mymodel1}
{\rm vec}(\bm{Y}_t) =(\bm{B}_1\otimes \bm{A}_1){ \circ}_{L} {\rm vec}(\bm{Y}_{t-1})+\cdots+(\bm{B}_p\otimes \bm{A}_p){ \circ}_{L} {\rm vec}(\bm{Y}_{t-p}) +{\rm vec}(\bm{\mathcal{E}}_{t}) ,~t \in \mathbb{Z},
\end{equation}
where ``$\otimes$'' is Kronecker product. 
If we further denote by $\bm{\Phi}_l=\bm{B}_l \otimes \bm{A}_l$ $(l\in\{1,...,p\})$,
then Model (\ref{Mymodel1}) can be  written as:
\begin{align}\label{vector_model}
{\rm vec}(\bm{Y}_t) =\bm\Phi_1{\circ}_{L} {\rm vec}(\bm{Y}_{t-1})+\cdots+ \bm\Phi_p{\circ}_{L} {\rm vec}(\bm{Y}_{t-p})+{\rm vec}(\bm{\mathcal{E}}_{t}) ,~t \in \mathbb{Z}.
\end{align}
It is clear that Model (\ref{vector_model}) takes the form as a multivariate general INAR (MGINAR($p$)) model of \cite{Latour1997}, which
builds a bridge  between the proposed MAT-INAR($p$) model and the MGINAR($p$) model.
However, Model (\ref{vector_model}) has  $p(mn)^2+nm$ parameters to be estimated,
which is a lot more than Model (\ref{Mymodelp}) who  contains only $p(m^2+n^2)+mn$ parameters.
Besides, Model (\ref{Mymodelp}) is more interpretable than the MGINAR($p$) model.
It  maintains the matrix structure of the observed data by the thinning operators defined in (\ref{circ}) and (\ref{circelement}).
Similar to MAR model (\cite{Chen2021}),
the coefficient matrices $\bm{A}_l$ and $\bm{B}_l$ of Model (\ref{Mymodelp}) are crucial.
The left matrix $\bm{A}_l$ reflects row-wise interactions, and the right matrix $\bm{B}_l$
introduces column-wise dependence.
Therefore,  Model (\ref{Mymodelp})
not only takes into account the interactive correlations between rows and columns of matrix-variate,
but also describes the autocorrelation of matrix time series across time.
Furthermore, we give some interpretations for Model (\ref{Mymodelp}) in Section C1 of Supplement Material.
If Model (\ref{vector_model}) is used for modeling,
the structure of the matrix data will be destroyed, resulting in serious information loss.

{\subsection{Probabilistic properties}\label{model_properties}

In the following proposition, we state the existence, strict stationarity and ergodicity of the MAT-INAR($p$) process defined in (\ref{Mymodelp}). 
To this end, we define  a $mnp \times mnp$ matrix ${\bm{\mathcal{A}}}$, where each element is a $mn \times mn$ block matrix  ${\bm{\mathcal{A}}}_{i,j}$ defined as  
\begin{align}\label{AAA}
&{\bm{\mathcal{A}}}_{i,j}=\bm{B}_j\otimes \bm{A}_j,~i=1,j\in\{1,...,p\},\nonumber\\
&{\bm{\mathcal{A}}}_{i,i-1}=\bm{I}_{mn},~i\in\{2,...,p\},~\text{otherwise},~{\bm{\mathcal{A}}}_{i,j}=\bm{0}_{mn},
\end{align}
where $\bm{I}_{mn}$ is a $(mn)$th-order identity matrix and
$\bm{0}_{mn}$ is a $(mn)$th-order zero matrix.
We will refer to $\bm{I}_{q\times q}$ as $\bm{I}_q$ unless there is any ambiguity.
Similarly,
we use $\bm{0}_q$ to denote a $q$th-order  square matrix with all elements equal to 0, instead of
using $\bm{0}_{q\times q}$ for simplicity. 
\begin{proposition}\label{stationary}
If   $\rho({\bm{\mathcal{A}}})<1$, then the MAT-INAR($p$) process  $(\ref{Mymodelp})$ is stationary and causal, where for any square matrix,  $\rho(\cdot)$   denotes its spectral radius, i.e., the maximum modulus of the (complex) eigenvalues of this matrix.
\end{proposition}

When the condition of Proposition \ref{stationary} is satisfied, 
we first derive the conditional expectation and expectation for the MAT-INAR($p$) process $(\ref{Mymodelp})$ as
\begin{flalign}
&\mbox{\textbullet ~} E(\bm{Y}_t|\bm{Y}_{t-1},...,\bm{Y}_{t-p}) = \sum_{l=1}^p\bm{A}_l \bm{Y}_{t-l}\bm{B}_l^{\top}+ \bm{\Lambda} ,~t \in \mathbb{Z}, \label{e1}&\\
&{\mbox{\textbullet ~}{\rm vec}(\bm{\mu}_Y):= E({\rm vec}(\bm{Y}_t)) =\left(\bm{I}_{mn}- \sum_{l=1}^p\bm{B}_l\otimes \bm{A}_l\right)^{-1} {\rm vec}(\bm{\Lambda}),~t \in \mathbb{Z}}.\label{e2}
\end{flalign}
Equation (\ref{e1}) can be obtained from   ({\bf P}B.4) in  Section B.1 of Supplement Material.
The conditional expectation combines the row-wise and column-wise interactions simultaneously.
Based on \cite{Latour1997}, equation (\ref{e2}) can be readily proven for Model (\ref{Mymodel1}).

Next, we go on to derive the variance-covariance matrix of $\bm{Y}_t$.
One purpose of introducing the MAT-INAR($p$) model is to be able to capture the cross-correlation of the
rows and columns of  MITS, so that the effect of all sequences can be considered simultaneously.
We  straighten the matrix-variate and define the vector version of the covariance as
$$
\bm{\Gamma}_0 := {\rm Cov}({\rm vec}(\bm{Y}_t),{\rm vec}(\bm{Y}_t))=E\left(({\rm vec}(\bm{Y}_t)-{\rm vec}(\bm{\mu}_Y))({\rm vec}(\bm{Y}_t)-{\rm vec}(\bm{\mu}_Y))^{\top}\right),~t\in \mathbb{Z}.
$$
However, this definition cannot well capture the correlation between the rows of $\bm{Y}_t$.
Thus,
we define the variance-covariance matrix of $\bm{Y}_t$ using the Kronecker product as
\begin{equation}\label{gamma01}
\bm{\Gamma}_0^{\otimes}=E\left((\bm{Y}_t-\bm{\mu}_Y)\otimes(\bm{Y}_t-\bm{\mu}_Y)^{\top}\right),~t \in \mathbb{Z}.
\end{equation}
Note that $\bm{Y}_t$ is a $m\times n$ matrix, then  $\bm{\Gamma}_0$ and $\bm{\Gamma}_0^{\otimes}$ are both $mn\times mn$ matrices.
The difference between $\bm{\Gamma}_0$ and $\bm{\Gamma}_0^{\otimes}$ is that the matrix elements are arranged differently.
Therefore, $\bm{\Gamma}_0^{\otimes}$ can not be used directly either.
To overcome this obstacle,
we introduce a transformation matrix $\bm{\mathcal{T}}=(t_{i,j})_{mn\times mn}$ (\cite{Samadi2014}), where $t_{i,j}$ is defined as $t_{i,j}=1$, if $i\in\{sm+1,sm+2,...,(s+1)m\}, ~j\in\{(s+1)+(i-sm-1)n\}$ and $s\in\{0,1,...,(n-1)\}$,  otherwise, $t_{i,j}=0$.
Based on the matrix $\bm{\mathcal{T}}$, we can transform  $\bm{\Gamma}_0^{\otimes}$ to several blocks such that
each block can be the variance of a vector, or the covariance of two rows or columns.
Then, we have the following column-wise or row-wise variance-covariance of $\bm{Y}_t$ defined as
\begin{equation}\label{leftcov}
\bm{\Sigma}_0^{c}:=E\left[\bm{\mathcal{T}}((\bm{Y}_t-\bm{\mu}_Y)\otimes(\bm{Y}_t-\bm{\mu}_Y)^{\top})\right]=\bm{\mathcal{T}}\bm{\Gamma}_0^{\otimes}
=(\bm{\Sigma}_{ij}^{c})_{n\times n},
\end{equation}
$\bm{\Sigma}_{i,j}^{c}$ is a $m\times m$ matrix, $\bm{\Sigma}_{i,j}^{c}=\bm{\Sigma}_{j,i}^{c}={\rm Cov}(\bm{Y}_{\cdot,i,t},\bm{Y}_{\cdot,j,t})$, $\bm{Y}_{\cdot,j,t}$ is the $j$th column of $\bm{Y}_t$, and
\begin{equation}\label{rightcov}
\bm{\Sigma}_0^{r}:=E\left[((\bm{Y}_t-\bm{\mu}_Y)\otimes(\bm{Y}_t-\bm{\mu}_Y)^{\top})\bm{\mathcal{T}}\right]=\bm{\Gamma}_0^{\otimes}\bm{\mathcal{T}}=(\bm{\Sigma}_{ij}^{r})_{m\times m},
\end{equation}
where $\bm{\Sigma}_{i,j}^{r}$ is a $n\times n$ matrix,
$\bm{\Sigma}_{i,j}^{r}=\bm{\Sigma}_{j,i}^{r}={\rm Cov}(\bm{Y}_{i,\cdot,t},\bm{Y}_{j,\cdot,t})$, $\bm{Y}_{i,\cdot,t}$ is the $i$th row of $\bm{Y}_t$.
We take the column-wise transformation  as an example.
As  seen in (\ref{leftcov}),  we transform $\bm{\Sigma}_0^{c}$  into $n^2$ sub-matrices, each of $m\times m$ dimension.
The diagonal block $\bm{\Sigma}_{j,j}^c$ is the variance matrix of column $\bm{Y}_{\cdot,j,t}$, while
the off-diagonal block $\bm{\Sigma}_{i,j}^c$ is the covariance matrix of  columns $\bm{Y}_{\cdot,i,t}$ and $\bm{Y}_{\cdot,j,t}~(i,j \in \{1,...,n\},~i\neq j)$.

Similarly, we can define the lag-$h$ autocovariance matrix of the MAT-INAR($p$) process with expression $\bm{\Gamma}_h^{\otimes}$ as follows $\bm{\Gamma}_h^{\otimes}=E\left((\bm{Y}_{t+h}-\bm{\mu}_Y)\otimes(\bm{Y}_t-\bm{\mu}_Y)^{\top}\right),~t \in \mathbb{Z}$.
Thus, transform $\bm{\Gamma}_h^{\otimes}$ by left multiplying and right multiplying $\bm{\mathcal{T}}$,
we obtain the lag-$h$ column-wise and row-wise autocovariance matrix in the forms of
$\bm{\Sigma}_h^{c}=\bm{\mathcal{T}}\bm{\Gamma}_h^{\otimes}$ and $\bm{\Sigma}_h^{r}=\bm{\Gamma}_h^{\otimes}\bm{\mathcal{T}}$, respectively.
It is worth noting that for a MITS $\{\bm{Y}_t\}$, $\bm{\Gamma}_h^{\otimes}\neq\bm{\Gamma}_{-h}^{\otimes}$ (\cite{Samadi2014}).
The relationship between $\bm{\Gamma}_h^{\otimes}$ and $\bm{\Gamma}_{-h}^{\otimes}$ is given by
${\bm{\Gamma}_h^{\otimes}}=(\bm{\mathcal{T}}\bm{\Gamma}_{-h}^{\otimes}\bm{\mathcal{T}})^{\top}$.

\section{Projection estimation}\label{sec:3}

\indent Let $\{\bm{Y}_t\}_{t=1}^T$ be a series of matrix-valued observations  generated from the  MAT-INAR($p$) process.
Denote by
$$
\bm{\Theta}=
\{
\underbrace{{[0,1]}^{m\times m}\times{[0,1]}^{n\times n}\times \cdots \times{[0,1]}^{m\times m}\times{[0,1]}^{n\times n}}_{p~\text{times the product of}~{[0,1]}^{m\times m}\times{[0,1]}^{n\times n}}\times(0,\infty)^{m\times n}
\}
$$
the
parameter space for parameters $\bm{A}_1,\bm{B}_1,...,\bm{A}_p, \bm{B}_p,\bm{\Lambda}$.
In the following, we study the projection (PROJ) estimation  and the asymptotic properties of the PROJ estimators $\widehat{\bm{A}}_{l,pr}$, $\widehat{\bm{B}}_{l,pr}~(l\in\{1,...,p\})$ and $\bm{\Lambda}_{pr}$.

\subsection{Projection method}
In this section, we consider the PROJ estimation for model parameters by two steps.
First, we obtain the least squares parameter estimates of the  vectorized MAT-INAR model, i.e, Model (5).
Second, based on the above least squares parameter estimates, we solve the nearest Kronecker product (NKP) problem to obtain $\widehat{\bm{A}}_{l,pr}$ and $\widehat{\bm{B}}_{l,pr}~(l\in\{1,...,p\})$.

In the first step, let recall that $\bm{\Phi}_l=\bm{B}_l \otimes \bm{A}_l~(l\in\{1,...,p\})$. We
derive the conditional least squares (CLS) estimates $\widehat{\bm{ \Phi}}_l$, $\widehat{\bm\Lambda}$ of $\bm{\Phi}_l$, $\bm{\Lambda}$ based on  Model  \eqref{vector_model}.
We denote  $\bm{\Psi}^{\top}:=(\bm\Phi_1,...,\bm\Phi_p, {\rm vec}(\bm\Lambda))$,
$\bm{X}_{t}^{\top}:=({\rm vec}(\bm{Y}_{t})^{\top},...,{\rm vec}(\bm{Y}_{t-p+1})^{\top}, {1})$,
${\bm{\mathcal X}}=(\bm{X}_{p},...,\bm{X}_{T-1})^{\top}$, and $\bm{\mathcal Y}=({\rm vec}(\bm{Y}_{p+1}),...,{\rm vec}(\bm{Y}_{T}))^{\top}$.
Then, the CLS criterion function of $\bm\Psi$ takes the form
$
Q(\bm\Psi):=\left(\bm{\mathcal Y}-G(\bm\Psi)\right)^{\top}\left(\bm{\mathcal Y}-G(\bm\Psi)\right),
$
where $G(\bm\Psi):=E(\bm{\mathcal Y}|{\rm vec}(\bm{Y}_{t-1}),...,{\rm vec}(\bm{Y}_{t-p}))={\bm{\mathcal X}}\bm\Psi$.
Then, the CLS-estimator $\widehat{\bm\Psi}$ can be obtained as $
\widehat{\bm{\Psi}}=\mathop{\arg\min}\limits_{\bm{\Psi} \in \bm{\Theta}^{\ast}}Q(\bm{\Psi}),
$
where $\bm{\Theta}^{\ast} \subseteq \mathbb{R}^{(pmn+1)\times (mn)}$ denotes the parameter space of $\bm{\Psi}$.
Solving the score equation ${\partial Q(\bm{\Psi})}/{\partial \bm{\Psi}}=\bm{0}$ gives
\begin{equation}\label{clspsi}
\widehat{\bm{\Psi}}=(\bm{\mathcal  X}^{\top}\bm{\mathcal X})^{-1}\bm{\mathcal  X}^{\top}\bm{\mathcal Y}.
\end{equation}
Thus, we obtain  $\widehat{\bm \Phi}_1,...,\widehat{\bm \Phi}_p$ and $\widehat{\bm\Lambda}$,
and ${\widehat{\bm\Lambda}}_{pr}$ is trivially obtained as $\widehat{\bm\Lambda}_{pr}=\widehat{\bm\Lambda}$.

In the second step, for  $\bm{A}_l$ and $\bm{B}_l$ ($l\in \{1,...,p\}$), we find the estimators  $\widehat{\bm{A}}_l$ and $\widehat{\bm B}_l$  by projecting $\widehat{\bm \Phi}_l$ onto the space of Kronecker products under Frobenius norms:
\begin{equation}\label{NKP}
(\widehat{\bm{A}}_l,\widehat{\bm{B}}_l)=\mathop{\arg\min}\limits_{\bm{A}_l, \bm{B}_l}{\left\Vert \widehat{\bm\Phi}_l-\bm{B}_l\otimes \bm{A}_l\right\Vert}_F^2,~l\in \{1,...,p\},
\end{equation}
which is known as the NKP problem (\cite{Loan2000}).
To obtain the solution of  (\ref{NKP}), 
we adopt the re-arrangement operator (\cite{Chen2021}) defined as
$
g:\mathbb{R}^{mn}\times \mathbb{R}^{mn}\rightarrow \mathbb{R}^{m^2}\times \mathbb{R}^{n^2}$. 
Detailed definition and related properties  about this operator can refer to Section B.3 of Supplement Material.

Based on the  properties  of the operator $
g$, the NKP problem (\ref{NKP}) is transformed to
\begin{align}
(\widehat{\bm{A}}_l,\widehat{\bm{B}}_l)=\mathop{\arg\min}\limits_{\bm{A}_l,\bm{B}_l}{\Vert \widehat{\bm\Phi}_l-\bm{B}_l\otimes\bm{A}_l\Vert}_F^2&=\mathop{\arg\min}\limits_{\bm{A}_l,\bm{B}_l}{\Vert g(\widehat{\bm\Phi}_l)-g(\bm{B}_l\otimes\bm{A}_l)\Vert}_F^2\notag\\
&=\mathop{\arg\min}\limits_{\bm{A}_l,\bm{B}_l}{\Vert \widetilde{\bm\Phi}_l-{\rm vec}(\bm{A}_l){\rm vec}(\bm{B}_l)^{\top}\Vert}_F^2,~l\in \{1,...,p\}.
\notag
\end{align}
Following \cite{Loan2000}, we can obtain its solution 
through the singular value decomposition (SVD) of $\widetilde{\bm\Phi}_l$,
which gives
$
{\rm vec}({\bm{{A}}}_{l})=\sqrt{d_{1}^{(l)}}\bm{u}_{1}^{(l)}$,
${\rm vec}({\bm{{B}}}_{l})^{\top}=\sqrt{d_{1}^{(l)}}{\bm{v}_{1}^{(l)}}^{\top}~(l\in\{1,...,p\}),
$
where $d_{1}^{(l)}$ is the largest singular value of $\widetilde{\bm{\Phi}}_l$, 
and $\bm{u}_{1}^{(l)}$ is the corresponding first left singular vector,  $\bm{v}_{1}^{(l)}$ is the corresponding first right singular vector. 
By the property of SVD, we have  ${\Vert\bm{{ A}}_{l}\Vert}_F=({d_{1}^{(l)}})^{1/2}$. 
Then, the estimators ${\widehat{\bm {A}}}_{l,pr}$ and ${\widehat{\bm {B}}}_{l,pr}$ are obtained by converting
$\bm{u}_{1}^{(l)}$ and
$d_1^{(l)}\bm{v}_1^{(l)}$  into matrices. 
Notice that ${\Vert\bm{{ A}}_{l}\Vert}_F=({d_{1}^{(l)}})^{1/2}$ implies ${\Vert\widehat{\bm{ A}}_{l,pr}\Vert}_F=1$ for $l\in\{1,...,p\}$, 
which ensures the identifiability of  Model {\rm (\ref{Mymodelp})}.

\subsection{The asymptotic property of PROJ estimator}

To present the asymptotic property of the projection estimators ${\widehat{\bm{A}}}_{l,pr}$ and ${\widehat{\bm {B}}}_{l,pr}$, we first give Theorem \ref{thm1} to guarantee $\widehat{\bm \Phi}_l$  converges to a multivariate normal distribution.

\begin{theorem}\label{thm1}
Let $\{\bm{Y}_t\}$ be a MAT-INAR($p$) process satisfying {$\rho({\bm{\mathcal{A}}})<1$}. The parameters $\bm{A}_l,~\bm{B}_l~(l\in\{1,...,p\})$, $\bm{\Lambda}$ and the covariance matrix of $\bm{\mathcal{E}}_{t}$ are nonsingular.
Then, 
the estimator $\widehat{\bm \Phi}_l$ implied by (\ref{clspsi}) is strongly consistent and asymptotically normal, i.e.,
\begin{equation}
\sqrt{T-p}~{\rm vec}(\widehat{\bm{\Phi}}_l^{\top}-\bm{\Phi}_l^{\top})\overset{L}{\longrightarrow}N(\bm{0},\bm{W}_l{\widetilde{\bm{\Sigma}}}\bm{W}_l^{\top}),~l\in\{1,...,p\},
\end{equation}
where $\bm{W}_l=(\bm{0}_{m^2n^2},...,\bm{I}_{m^2n^2},...,\bm{0}_{m^2n^2},\bm{0}_{mn})$, the $l$th block element is
a $(mn)^2$th-order identity matrix,
and the rest is a zero matrix, $\widetilde{\bm\Sigma}= \bm{\Sigma}_{\mathcal  U}\otimes\bm{H}^{-1} $,
$\bm{\Sigma}_{\mathcal{U}}:=E(\bm{\mathcal U}_{t}\bm{\mathcal U}_{t}^{\top})$ with
$\bm{{\mathcal{U}}}_t={\rm vec}(\bm{Y}_t)-{\rm vec}(\bm{{\Lambda}})-\sum_{l=1}^p{\bm {\Phi}}_l{\rm vec}(\bm{Y}_{t-l})
$,
$\bm{H}:=E(\bm{ X}_{t}\bm{ X}_{t}^{\top})$.
Moreover, $\widehat{\bm\Sigma}_{\mathcal{U}}=\sum_{t=p+1}^T\widehat{\bm {\mathcal{U}}}_t\widehat{\bm {\mathcal{U}}}_t^{\top}/(T-mn-p)$ converges a.s. to $\bm{\Sigma}_{\mathcal{U}}$,
and $\widehat{\bm{ H}}_T=\sum_{t=p+1}^T\bm{X}_{t}\bm{X}_{t}^{\top}/(T-p)$ converges a.s. to $\bm{H}$.
\end{theorem}

Following the standard theory of MGINAR($p$) model in \cite{Latour1997},
we can easily prove that Theorem \ref{thm1} holds.
Now, we study the asymptotic property of  the PROJ-estimators.
For this purpose, we need to introduce some notations.
For $l\in\{1,...,p\}$, we define $\bm{\alpha}_l:={\rm vec}(\bm{A}_l)$, 
$\bm{\beta}_l:={\rm vec}(\bm{B}_l)$, 
$\bm{\beta}_l^{(1)}:=\bm{\beta}_l/ {\Vert\bm{\beta}_l\Vert}$. Thus,
we have that $\bm{\alpha}_l$ and $\bm{\beta}_l^{(1)}$ are unit vectors.
Then the following results establish the strong consistency and the asymptotic normality of the PROJ-estimators.

\begin{theorem}\label{thm2}

Under  the conditions of Theorem \ref{thm1},
the PROJ-estimators $\widehat{\bm{A}}_{l,pr}$ and $\widehat{\bm{B}}_{l,pr}$  are strongly consistent and asymptotically normal, i.e.,
\begin{equation}\label{thm2.1}
\sqrt{T-p}\left(\begin{array}{c}
{\rm vec}(\widehat{\bm{A}}_{l,pr}-\bm{A}_l) \\
{\rm vec}(\widehat{\bm{B}}_{l,pr}-\bm{B}_l)
\end{array}\right)\overset{L}{\longrightarrow}N\left(\bm{0},\bm{V}_0^{(l)}\bm{\Xi}_1^{(l)}{\bm{V}_0^{(l)}}^{\top}\right),~l\in\{1,...,p\},
\end{equation}
where 
\begin{align}
&\bm{V}_{0}^{(l)}:=\left(\begin{array}{c}
{\Vert\bm{B}_l\Vert}_F^{-1} [{\bm{\beta}_l^{(1)}}^{\top}\otimes(\bm{I}_{m^2}-\bm{\alpha}_l\bm{\alpha}_l^{\top})]\nonumber\\
\bm{I}_{n^2}\otimes\bm{\alpha}_l^{\top}
\end{array}\right),\nonumber
\end{align}
{$\bm{\Xi}_1^{(l)}=\widetilde{\bm{\mathcal{T}}}\bm{W}_l\widetilde{\bm\Sigma}^{\top}\bm{W}_l^{\top}\widetilde{\bm{\mathcal{T}}}$,  $\widetilde{\bm{\mathcal{T}}}$ is defined in Section D3  of Supplement Material.
Moreover, 
\begin{equation}\label{thm2.2}
\sqrt{T-p}\left({\rm vec}(\widehat{\bm{B}}_{l,pr})\otimes{\rm vec}(\widehat{\bm{A}}_{l,pr})-{\rm vec}(\bm{B}_l)\otimes{\rm vec}(\bm{A}_l)\right)\overset{L}{\longrightarrow}N\left(\bm{0},\bm{V}_1^{(l)}\bm{\Xi}_1^{(l)}{\bm{V}_1^{(l)}}^{\top}\right),
\end{equation}
where
$\bm{V}_1^{(l)}:=(\bm{\beta}_l^{(1)}{\bm{\beta}_l^{(1)}}^{\top})\otimes\bm{I}_{m^2}+\bm{I}_{n^2}\otimes(\bm{\alpha}_l\bm{\alpha}_l^{\top})-(\bm{\beta}_l^{(1)}{\bm{\beta}_l^{(1)}}^{\top})\otimes(\bm{\alpha}_l\bm{\alpha}_l^{\top})$.
}
\end{theorem}

\begin{theorem}\label{thm3}
Under  the conditions of Theorem \ref{thm1}, the PROJ-estimator $\widehat{\bm \Lambda}_{pr}$ is strongly consistent and asymptotically normal, i.e.,
$
\sqrt{T-p}~{\rm vec}(\widehat{\bm \Lambda}_{pr}^{\top}-\bm{ \Lambda}^{\top})\overset{L}{\longrightarrow}N\left(\bm{0},\bm{W}_{p+1}\widetilde{\bm\Sigma}\bm{W}_{p+1}^{\top}\right),
$
where $\bm{W}_{p+1}=(\bm{0}_{m^2n^2},...,\bm{0}_{m^2n^2},\bm{I}_{mn})$.
\end{theorem}

The  vectorization structure (\ref{vector_model}) is useful in deriving the initial estimators and in study the asymptotic properties of the projection estimators. 
However, the projection estimator as well as Theorems \ref{thm2} and  \ref{thm3} still require that the observed matrix-valued time series satisfying Model (\ref{Mymodelp}).  
On the premise that Theorem \ref{thm1} is true, we can easily prove Theorem \ref{thm2} by extending the arguments in the proof of  Theorem 2 in \cite{Chen2021}. 

\section{Iterated conditional least squares estimation}\label{sec:4}
When  the dimension $(pmn+1)\times(mn) $ of $\bm{\Psi}$ is high,
the problem will become a high-dimensional or even ultra-high-dimensional estimation problem.
In this situation, the resulting projection estimators may not be very accurate.
In this section,
we derive the iterated conditional least squares (ICLS) estimation for parameters $\bm{A}_l,\bm{B}_l~(l\in\{1,...,p\})$ and $\bm{\Lambda}$.
In this procedure, the PROJ estimators are taken as the initial values,
which can help to find a more accurate estimator based on the iterative algorithm.

The ICLS estimators $\widehat{\bm{A}}_{l,ils}$, $\widehat{\bm{B}}_{l,ils}~(l\in\{1,...,p\})$ and $\widehat{\bm{\Lambda}}_{ils}$ are the
solution of
\begin{equation}\label{ICLS}
(\widehat{\bm{A}}_1, \widehat{\bm{B}}_1,...,\widehat{\bm{A}}_p,\widehat{\bm{B}}_p,\widehat{\bm{\Lambda}})=
\mathop{\arg\min}\limits_{\bm{A}_1,\bm{B}_1,...,\bm{A}_p,\bm{B}_p,\bm{\Lambda}}\sum_{t=p+1}^{T}\left\Vert \bm{Y}_t-\sum_{l=1}^p\bm{A}_l\bm{Y}_{t-l}\bm{B}_l^{\top}-\bm{\Lambda}\right\Vert_F^2.
\end{equation}
Taking partial derivatives of (\ref{ICLS}) with respect to  $\bm{A}_l$, $\bm{B}_l$ ($l\in\{1,...,p\}$) and $\bm{\Lambda}$,
we obtain the gradient conditional for the ICLS as follows:
\begin{equation}\label{iteration}
    \begin{cases}
    \sum_{t=p+1}^T \left(\bm{Y}_{t}- \sum_{l=1}^p\bm{A}_l\bm{Y}_{t-l}{\bm{B}_l}^{\top}- \bm{\Lambda}\right)\bm{B}_k\bm{Y}_{t-k}^{\top}=\bm{0},~k\in\{1,...,l\},\\
    \sum_{t=p+1}^T \left(\bm{Y}_{t}^{\top}- \sum_{l=1}^p{\bm{B}_l}\bm{Y}_{t-l}^{\top}\bm{A}_l^{\top}-\bm{\Lambda}^{\top}\right)\bm{A}_k\bm{Y}_{t-k}=\bm{0},~k\in\{1,...,l\},\\
   (T-p)\bm{\Lambda}- \sum_{t=p+1}^T  \left(\bm{Y}_{t}- \sum_{l=1}^p\bm{A}_l\bm{Y}_{t-l}{\bm{B}_l}^{\top}\right)=\bm{0}.
    \end{cases}
\end{equation}
Considering the complex product structures of the parameter matrices in (\ref{iteration}), it is not easy to find the closed-form solution.
Therefore, we iteratively update one, while keeping the other fixed, these iterations are given by following three steps:

{\bf Step 1}: Given initial values $\bm{A}_l^{(0)}=\widehat{\bm{A}}_{l,pr}$, $\bm{B}_l^{(0)}=\widehat{\bm{B}}_{l,pr}~(l\in\{1,...,p\})$ and
$\bm{\Lambda}^{(0)}=\widehat{\bm{\Lambda}}_{pr}$. 

{\bf Step 2}: In the $s$th step, update $\bm{A}_l^{(s)}$, $\bm{B}_l^{(s)}$ ($l\in\{1,...,p\}$), $\bm{\Lambda}^{(s)}$  according to (\ref{iteration}) to get
\begin{align}
\bm{A}_l^{(s+1)}&\leftarrow\bm{M}_l^{(s)}{\left(\sum_{t=p+1}^T\bm{Y}_{t-l}{ \bm{B}_l^{(s)}}^{\top} \bm{B}_l^{(s)}  \bm{Y}_{t-l}^{\top}\right)^{-1}},\nonumber\\
\bm{B}_l^{(s+1)}&\leftarrow\bm{N}_l^{(s)}{\left(\sum_{t=p+1}^T\bm{Y}_{t-l}^\top{ \bm{A}_l^{(s+1)}}^{\top} \bm{A}_l^{(s+1)}  \bm{Y}_{t-l}\right)^{-1}},\nonumber\\
\bm{\Lambda}^{(s+1)}&\leftarrow \frac{1}{(T-p)}\sum_{t=p+1}^T\left(\bm{Y}_{t}-\sum_{l=1}^p{\bm{A}_l^{(s+1)}}\bm{Y}_{t-l} {\bm{B}_l^{(s+1)}}^{\top}\right),\nonumber
\end{align}
where
$$
\bm{M}_l^{(s)}:=\sum_{t=p+1}^T \left(\bm{Y}_{t} -\sum_{k=1}^{l-1}\bm{A}_k^{(s+1)}\bm{Y}_{t-k}{\bm{B}_k^{(s)}}^{\top}-\sum_{k=l+1}^{p}\bm{A}_k^{(s)}\bm{Y}_{t-k}{\bm{B}_k^{(s)}}^{\top}- \bm{\Lambda}^{(s)}\right){\bm{B}_l^{(s)} \bm{Y}_{t-l}^{\top}},
$$
$$
\bm{N}_l^{(s)}:=\sum_{t=p+1}^T \left(\bm{Y}_{t}^\top -\sum_{k=1}^{l-1}{\bm{B}_k^{(s+1)}}^{\top}\bm{Y}_{t-k}^\top\bm{A}_k^{(s+1)}-\sum_{k=l+1}^{p}{\bm{B}_k^{(s)}}^{\top}\bm{Y}_{t-k}^\top\bm{A}_k^{(s+1)}- {\bm{\Lambda}^{(s)}}^\top\right){\bm{A}_l^{(s+1)} \bm{Y}_{t-l}}.
$$

{\bf Step 3}: Repeat step 2  until
$
\max\left\{\left\| \bm{A}_l^{(s+1)}-\bm{A}_l^{(s)}\right\|_F, \left\| \bm{B}_l^{(s+1)}-\bm{B}_l^{(s)}\right\|_F,
\left\| \bm{\Lambda}^{(s+1)}-\bm{\Lambda}^{(s)}\right\|_F\right\} <c\times {10}^{-\delta},
$
for some positive constants $c$ and $\delta$. Thus, we obtain  $\widehat{\bm {A}}_{l,ils}$, $\widehat{\bm{B}}_{l,ils} ~(l\in\{1,...,p\})$ and $\widehat{\bm{\Lambda}}_{ils}$. In this study, we choose $c=1$ and $\delta=9$ without loss of generality.

The following theorem establishes the strong consistency and the asymptotic normality of the
ICLS estimators.

\begin{theorem}\label{thm4}
Let $\{\bm{Y}_t\}$ be a MAT-INAR($p$) process satisfying the conditions of Theorem \ref{thm1}, then we have
\begin{align}\label{thm3e}
&\sqrt{T-p}\left(
{\rm vec}(\widehat{\bm{A}}_{1,ils}-\bm{A}_1)^{\top},
{\rm vec}(\widehat{\bm{B}}_{1,ils}-\bm{B}_1),
\cdots,
{\rm vec}(\widehat{\bm{A}}_{p,ils}-\bm{A}_p)^{\top}, 
{\rm vec}(\widehat{\bm{B}}_{p,ils}-\bm{B}_p),\nonumber\right.\\
&\left.{\rm vec}(\widehat{\bm{V}}_{ils}-\bm{V})^{\top}
\right)^\top
\overset{L}{\longrightarrow}N(\bm{0},\bm{\Xi}_2),
\end{align}
where $\bm{\Xi}_2:=\bm{Q}^{-1}E(\bm{P}_t\bm{\Sigma}_{\mathcal U}\bm{P}_t^{\top})\bm{Q}^{-1}$, 
$\bm{Q}:=E(\bm{P}_t\bm{P}_t^{\top})+\sum_{l=1}^p\bm{\gamma}_l\bm{\gamma}_l^{\top}$, 
$$
\bm{P}_t:=\left(\bm{B}_1\bm{Y}_{t-1}^{\top}\otimes{\bm I}_m,{\bm I}_n\otimes\bm{A}_1\bm{Y}_{t-1},\cdots,\bm{B}_p\bm{Y}_{t-p}^{\top}\otimes{\bm I}_m,{\bm I}_n\otimes\bm{A}_p\bm{Y}_{t-p},{\bm I}_{mn}\right)^{\top},
$$  
and 
$
\bm{\gamma}_l=(\bm{0}_{1\times m^2 },\bm{0}_{1\times n^2 },...,\underbrace{\bm{\alpha}_{l}^{\top},\bm{0}_{1\times n^2}}_{l\text{th position} },...,\bm{0}_{1\times m^2},\bm{0}_{1\times n^2},\bm{0}_{1\times mn})^{\top}. 
$
\end{theorem}


\section{Determining the order for MAT-INAR($p$) model}\label{sec:5}
While the general MAT-INAR($p$) model provides
more flexibility and capability to capture different interactions among column and row variables, it also
poses the challenge of finding the order. 
This issue is crucial because accurate model order determination can effectively avoid problems such as overfitting or underfitting of the model. In this study, 
we propose an information criterion based procedure, which achieves selection consistency under some fixed initial order.

For a sufficiently large initial order ${\tilde p}$, we define the information criterion as
\begin{equation}\label{IC1}
IC_1(\tilde p)=\log\left(\frac{1}{T}\sum_{t={\tilde p}+1}^T{\left\Vert\bm{Y}_t-\sum_{l=1}^{\tilde p}\widehat{\bm{A}}_l\bm{Y}_{t-l}\widehat{\bm{B}}_l^{\top}-\widehat{\bm{\Lambda}}\right\Vert}_F^2\right)+\frac{1}{T}{\tilde p}\log T,
\end{equation}
where 
$T$ is the sample size, 
$\widehat{\bm{A}}_l$, $\widehat{\bm{B}}_l~(l\in\{1,...,{\tilde p}\})$ and $\widehat{\bm{\Lambda}}$ are the estimators obtained under the order ${\tilde p}$. 
In practice we typically cap the suitable order at some given $\overline{p}$, $\overline{p}$ is a given up bound to choose order. 
Then the estimated $\hat{p}$ is given by $\hat{p}={\arg\min}_{1\leq \tilde p \leq \overline{p}}IC_1(\tilde p)$.

Now, we state the consistency of order determination estimator in the following theorem. 
\begin{theorem}\label{thm5}
Let $\{\bm{Y}_t\}$ be a MAT-INAR($p$) process satisfying {$\rho({\bm{\mathcal{A}}})<1$}. 
Assuming that there exists some constant $\eta>0$ such that ${\Vert \bm{B}_l\otimes\bm{A}_l\Vert}_F^2\geq \eta$ for all $l\in\{1,...,p\}$ and ${\Vert\bm{\Lambda}\Vert}_F^2\geq \eta$, we have
$
\mathop{\lim}\limits_{T\rightarrow \infty}P(\hat{p}=p)=1.
$
\end{theorem}

The  criterion can be viewed as an extended Bayesian information criterion (\cite{Chen2008,Guo2016}). 
The difference between this criterion and the traditional BIC criterion is that when the model cannot calculate the likelihood function, the goodness of fit of the model is described by the residual sum of squares. 
With reference to this idea, the criterion can also be used for order selection by vector or univariate  model.
We present an alternative information  criterion which is  applicable to vectorized models in Section E of Supplement Material.

\section{Simulation studies}\label{sec:6}
In this section, we study the empirical performances of the proposed estimators and the order determination procedure.
The simulation studies are grouped into two part: first on the estimation accuracy, and
later on the order determination procedure.

\subsection{Comparisons of PROJ and ICLS}\label{sec:8.1}

\indent To report the performances of the proposed PROJ and ICLS methods mentioned previously, we conducted simulation
studies using $\textsf{R}$ software based on 1000 replications.
We consider the following  scenario 
under sample sizes $T = 200, 600, 1000$, respectively.

\begin{enumerate}[\text{Scenario} A.]
  \begin{sloppypar}
  \setlength{\itemsep}{1pt}
  \setlength{\parskip}{0pt}
  \setlength{\parsep}{0pt}
  \item We consider a MAT-INAR(1) model with $(m,n)=(2,2)$, and the initial parameters in $\bm{A}$, $\bm{B}$ and $\bm{\Lambda}$ are chosen as
\begin{equation}
 \widetilde{\bm{A}}= \left(\begin{array}{cc}
0.20 &0.40 \\
0.40& 0.20\\
\end{array}\right),~
\bm{B}= \left(\begin{array}{cc}
0.50 &0.30 \\
0.30& 0.50\\
\end{array}\right),~ \text{and}~
\bm{\Lambda}= \left(\begin{array}{cc}
1.00 &1.00\\
1.00&1.00\\
\end{array}\right), 
\notag
\end{equation}
and $\bm{A}=\widetilde{\bm{A}}/\Vert\widetilde{\bm{A}}\Vert_F$ to guarantee the uniqueness holds for the model. 
$\bm{\mathcal{E}}_{t}$ follows a matrix-variate Poisson  (Mpois) distribution (\cite{Yurchenko2021}) with mean $\bm\Lambda$, i.e., $\bm{\mathcal{E}}_{t}\sim Mpois(\bm\Lambda)$ given in Section B4 of Supplement Material.
\end{sloppypar}
\end{enumerate}
Scenarios B--C can be found at Section F in Supplement Material.

\begin{figure}[h]
\centering
\includegraphics[width=4.2in]{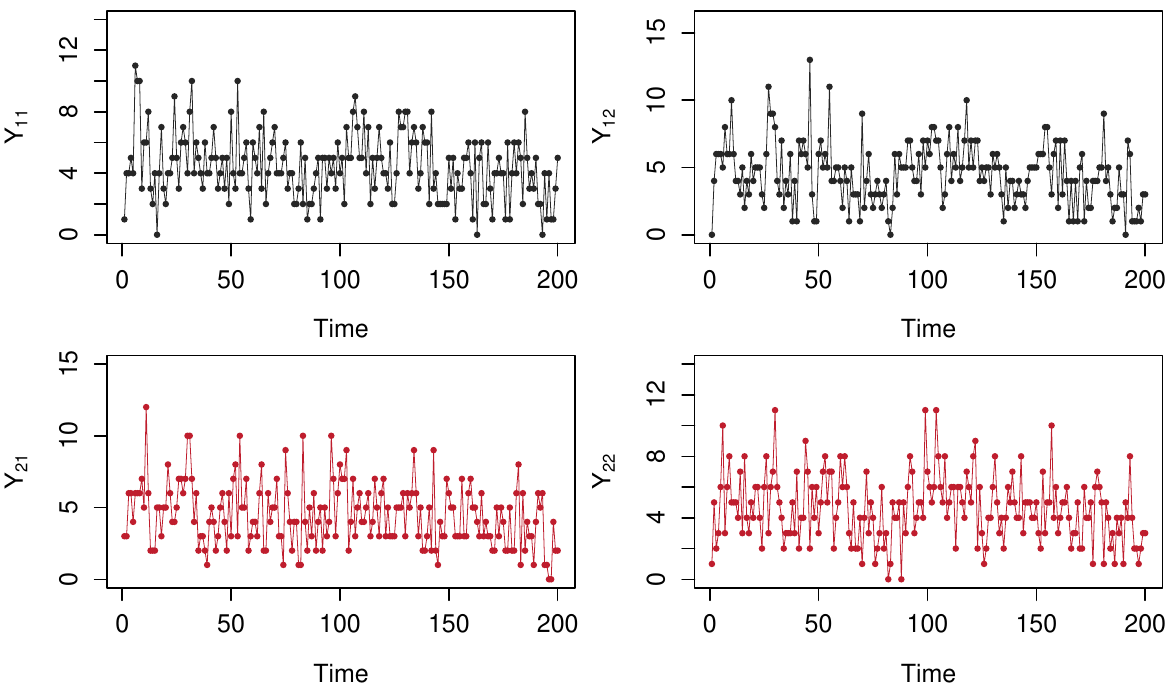}
\caption{
Time series plots of Scenario A.
}
\label{pathA}
\end{figure}

\begin{figure}[ht]
\centering
\includegraphics[width=6.in]{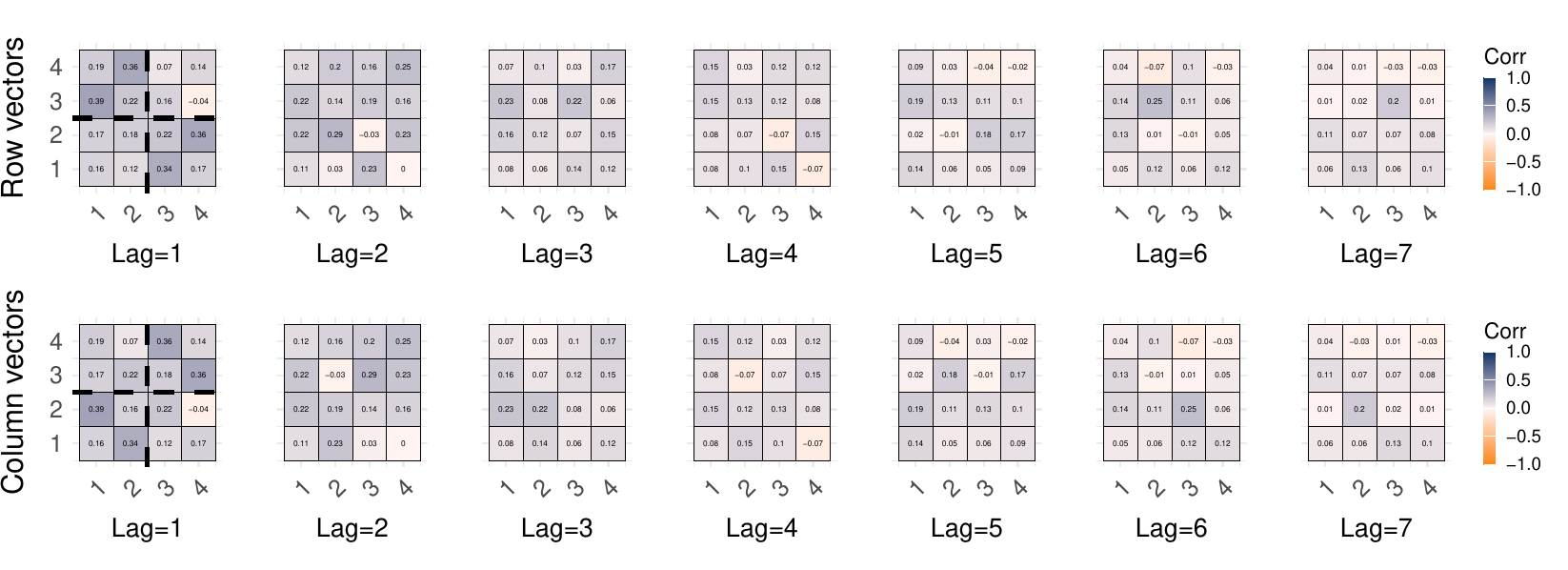}
\caption{
Cross-ACFs of row vectors and of column vectors  for Scenario A.
The black bold dotted line shows the ACF matrix blocks.
}
\label{aACF}
\end{figure}
Meanwhile, all elements of $\bm{A}_l$ and $\bm{B}_l$ are chosen to satisfy 
$\rho({\bm{\mathcal{A}}})<1$ to guarantee the fulfillment of the stationary condition. 
In order to investigate sample properties of  Scenario A, we draw the time series  in Figure \ref{pathA} and
cross-autocorrelation function (cross-ACF) plots  in Figure \ref{aACF} when $T=200$. 
We can see in Figure \ref{pathA} that each series fluctuates up and down around a fixed value within a constant range,
indicating that the generated series are stationary. 
Moreover, we find that in these simulated series, there  can
contain overdispersion, underdispersion and equidispersion components together, showing 
great flexibility for matrix-variate time series.   
Now we focus on the cross-ACFs given in Figure \ref{aACF}. 
The first panel on top is the cross-ACFs of the row vector of $\bm{Y}_t$, while the second panel on bottom denotes the cross-ACFs of the column vector of $\bm{Y}_t$. 
As  seen in Figure \ref{aACF}, (i) with the  increases of lag $h$, 
the color of the matrix block in each panel gradually becomes lighter, indicating that the autocorrelations of the matrix series are gradually weakening; 
(ii) none of the block shows a zero number, implying that the row and column sequences of the MAT-INAR process  has 
autocorrelation and interactions between each other. 

In order to investigate the performances of the proposed estimators, we calculate the empirical biases (Bias),  standard deviations (SD) of the estimates across 1000 replications, as well as the standard errors (SE) of the PROJ-estimators and the CLS-estimators. 
The simulation results of Scenario A are summarized in Table \ref{t1}. 

\begin{landscape}
\begin{table}[tbp]                    
\renewcommand\arraystretch{0.8}                 \setlength{\abovecaptionskip}{2pt}
\setlength{\belowcaptionskip}{10pt}
\caption{Simulation results for Scenario A: Bias, SE and SD}                 
\centering                                  
  \label{t1}
{\tabcolsep0.05in                           
\begin{tabular}{rrrrrrrrrrrrrrr}
  \hline
  Method& $T$  & Result & $a_{1,1}$ & $a_{2,1}$ &$a_{1,2}$ & $a_{2,2}$& $b_{1,1}$ & $b_{2,1}$ &$b_{1,2}$ & $b_{2,2}$& $\lambda_{1,1}$ & $\lambda_{2,1}$ &$\lambda_{1,2}$ & $\lambda_{2,2}$\\
  \hline
PROJ &200 &Bias&$-$0.013 &0.001 &$-$0.004 &$-$0.012 &$-$0.007 &$-$0.004 &0.000 &$-$0.002 &0.073 &0.056 & 0.070 &0.046\\
    &     &SD &0.076 &0.062 &0.062 &0.076 &0.064 &0.064 &0.061 &0.063 &0.365 &0.373 &0.371 &0.391\\
    &     &SE &0.074 &0.058 &0.058 &0.074 &0.062 &0.062 &0.062 &0.062 &0.407 &0.408 &0.408 &0.408\\
    &600  &Bias&$-$0.005 &0.001 &$-$0.001 &$-$0.005 &$-$0.003 &0.000 &0.001 &$-$0.002 &0.034 &0.011 &0.014 &0.019  \\
    &     & SD  & 0.043 & 0.036& 0.036& 0.044& 0.036 & 0.038 & 0.036& 0.036& 0.221& 0.226& 0.222& 0.233\\
    &     & SE  & 0.042 & 0.033& 0.033 & 0.042& 0.035 & 0.035 & 0.035& 0.035& 0.234 & 0.234 & 0.234& 0.233\\
    &1000  & Bias &$-$0.001 &0.000 & $-$0.002&$-$0.002 &$-$0.002 &0.001 & 0.000&0.000 &0.013 &0.005& 0.000&0.012  \\
    &     & SD  & 0.033& 0.026 & 0.027& 0.033& 0.028 & 0.028 & 0.028& 0.027& 0.173& 0.168 & 0.172& 0.178\\
    &     & SE  & 0.032 & 0.025& 0.025& 0.032& 0.027 & 0.027 & 0.027& 0.027& 0.181 & 0.181 & 0.180& 0.181\\
ICLS &200 &Bias&$-$0.012 &$-$0.000 &$-$0.003 &$-$0.013 &$-$0.009 &0.000 &0.002 &$-$0.010 &0.076 &0.066 &0.080 &0.070\\
    &     &SD &0.077 &0.075 &0.075 &0.073 &0.056 &0.059 &0.060 &0.058 &0.334 &0.332 &0.338 &0.322 \\
    &     &SE &0.074 &0.059 &0.059 &0.074 &0.063 &0.062 &0.062 &0.062 &0.333 &0.333 &0.331 &0.333 \\
    &600  & Bias &$-$0.005 &$-$0.003 & 0.000&$-$0.004 &$-$0.004&0.002 &0.000&$-$0.002&0.030 &0.033 & 0.021&0.022\\
    &     & SD  & 0.044 & 0.044 & 0.045 & 0.044 & 0.033 & 0.035 & 0.034 & 0.034 & 0.195 & 0.204 & 0.200 & 0.198\\
    &     & SE  & 0.043 & 0.034 & 0.034 & 0.043& 0.036 & 0.036& 0.036 & 0.036& 0.196 & 0.196 & 0.197 & 0.197\\
    &1000  & Bias &$-$0.002 &$-$0.002 & 0.000&$-$0.002&$-$0.001 &0.000 &0.000&$-$0.001 &0.008 &0.014 & 0.009&0.015\\
    &     & SD  & 0.032 & 0.034& 0.034 & 0.033 & 0.026 & 0.026& 0.027 & 0.026 & 0.153 & 0.151& 0.153 & 0.150 \\
    &     & SE  & 0.033 & 0.026 & 0.026 & 0.033 & 0.028 & 0.028& 0.028 & 0.028 & 0.153 & 0.153& 0.153& 0.153 \\
  \hline
\end{tabular}
}
\end{table}
\end{landscape}

As seen in Table \ref{t1}, all the Biases, SDs and SEs decrease as sample size $T$ increases.
This implies our estimators are consistent for all  parameters. 
We also see that the values of SD do not differ much from their corresponding SE, indicating that the 
estimators converge fast. 
In addition, most biases of the  ICLS-estimators  are smaller than the corresponding PROJ-estimators, and  the difference between the values of SD and their corresponding SE of ICLS-estimators  are smaller than the corresponding PROJ-estimators.
This means the ICLS-estimators perform better than the PROJ-estimators. 

\subsection{Autoregressive order $p$ estimation}\label{sec:8.3}

In this section, we conduct simulations to 
report the performances of the  information criteria (\ref{IC1}) introduced previously. 
Specifically, we consider three new scenarios, denoted by Scenarios D, E and F. 
For each scenario, the true orders of the model are $p =1,2,3,4$, respectively. 
For generating the data, 
the true values of parameters $\bm{A}_l$ and $\bm{B}_l~(l\in\{1,2,3,4\})$ in each scenario are generated randomly from  $U(0,1)$ distribution. 
Also, standardization for  $\bm{A}_l~(l\in\{1,2,3,4\})$ are applied to ensure $\Vert\bm{A}_l\Vert_F=1$. 
Moreover, the distribution settings and the parameter values of 
$\bm{\mathcal{E}}_{t}$ are chosen to be the same as Scenarios A, B and C. 
The simulation results are summarized in Tables \ref{t2}. 

\begin{table}[ht]
\renewcommand\arraystretch{0.7}      
\setlength{\belowcaptionskip}{0.2cm}
\caption{Relative frequencies of occurrence  of the events
$\{\hat{p}=p\}$,$~\{\hat{p} > p\}$ and $\{\hat{p} < p\}$ for IC$_1$.}
\label{t2}
\centering
{\tabcolsep0.03in             
\begin{tabular}{@{}ccccccccccccc@{}}
\hline
      &   & \multicolumn{3}{c}{Scenario D}      & &\multicolumn{3}{c}{Scenario E}& &\multicolumn{3}{c}{Scenario F}\\
      \cline{3-5}\cline{7-9}\cline{11-13}
    $T$  &   & $\{\hat{p}=p\}$ &$\{\hat{p}> p \}$ & $\{\hat{p}<p\}$ & &$\{\hat{p}=p\}$ &$\{\hat{p}> p\}$ & $\{\hat{p}<p\}$& &$\{\hat{p}=p\}$ &$\{\hat{p}> p\}$ & $\{\hat{p}<p\}$\\
  \hline
  200      & $p=1$ &  0.800    & 0.200    & 0.000  &&  0.782    & 0.218   & 0.000 &&  0.818    & 0.182   & 0.000 \\
          & $p=2$ & 0.583    & 0.118    & 0.299  && 0.522    &0.419   & 0.059&&  0.469    & 0.136   & 0.395  \\
          & $p=3$ &  0.506    & 0.105    &0.389  &&0.440   & 0.355    & 0.205 &&  0.535    & 0.087   & 0.378 \\
          & $p=4$ & 0.358   &0.039  & 0.603 && 0.533    & 0.174    & 0.293&&  0.473    & 0.165   & 0.362 \\
  600     & $p=1$ & 0.928    & 0.072  & 0.000  &&  0.959   & 0.041   & 0.000  &&  0.758    & 0.242   & 0.000\\
          & $p=2$ &0.939    &0.050   & 0.011  &&  0.923   & 0.031   & 0.046  &&  0.807    & 0.193   & 0.028\\
          & $p=3$ &0.941    &0.030    &0.029  &&  0.759    & 0.013   & 0.228 &&  0.891    & 0.107   & 0.002 \\
          & $p=4$ &0.810    &0.014    &0.176 &&  0.766   & 0.003    & 0.231 &&  0.879   & 0.037   & 0.084\\
  1000     & $p=1$ &0.959    &0.041   & 0.000  &&  0.988   &0.012 & 0.000&&  0.826   &0.174  & 0.000\\
          & $p=2$ &0.968    &0.031  & 0.001  &&0.992   &0.006  & 0.002&&  0.869   &0.131  & 0.000 \\
          & $p=3$ &0.982   &0.018  &0.000  &&  0.964    & 0.006    & 0.030&&  0.888    &0.112 & 0.000 \\
          & $p=4$ &0.971    &0.006   &0.023 && 0.964    &0.002    & 0.034 &&  0.969   &0.027  & 0.004\\
\hline
\end{tabular}
}
\end{table}

Table \ref{t2} reports the relative frequencies of occurrence of the events $\{\hat{p}=p\}$,$~\{\hat{p}> p\}$ and $\{\hat{p}<p\}$ 
with $\overline{p}=6$ for the information criterion (\ref{IC1}) under sample sizes $T=200,600$ and $1000$, 
while 
Table F5 in Section F shows the corresponding results for  the  criterion (E1) defined in Section E of Supplement Material.  
We can see  in Tables \ref{t2} and F5 that when the sample size is small ($T=200$), 
$\hat{p}$ estimated by (\ref{IC1}) chooses the true value $p$ at most times under different scenarios, 
while the estimated frequencies of $\hat{p}$ estimated by  the  criterion (E1) perform less well (see Table F5}), especially when 
$p=4$. 
This implies the  criterion (\ref{IC1})  outperforms the criterion (E1). 
However, as the sample size  $T$ increases, the probabilities that the order is correctly estimated gradually approach to one, 
namely, the two estimators $\hat{p}$ converge to true order $p$ as  $T$ increases. 
These results imply that the estimators of  the  criteria (\ref{IC1}) and  (E1) both  perform well when $T$ is large. 
When $T$ is small, we suggest to consider the criterion (\ref{IC1}) to obtain a more precise result for matrix-valued data.

\section{Application: Offences data in NSW}\label{sec:7}
\indent According to 2021 data released by World Health Organization (WHO)  and its partners, violence against women remains extremely prevalent and begins at a very young age.
In fact, there are many countries and regions have not done enough to protect women from domestic and sexual violence. Motivated by this fact, we consider a set of offence counts in Australia in this study.
Specifically,
we mainly consider data on two major crime types: sexual assault and domestic violence, 
aiming to study the patterns and intrinsic connections of such criminal activities. 

Many studies have shown that there are strong correlations between criminal activities in neighboring cities (\cite{Yang2023,Yang2024CXM}). 
Similar criminal acts and trans-regional criminal networks may lead to some criminal acts in neighboring regions.
Therefore, it is necessary to consider the joint modeling of these two types of crime in multiple cities. 
We select two adjacent cities in New South Wales (NSW), Australia: Waverley and Ryde.
As  seen in Figure G1 of Supplement Material, the two areas are located in the east coast of NSW.
Ryde is on the south side of the Parramatta River,
while Waverley is on the north side.

We choose the monthly  domestic violence assault and sexual offence counts of these two cities,
starting from  August 1996 to  August 2019, totally $T=277$ matrix-valued observations.
Without loss of generality, we denote the domestic violence assault counts in Waverley by $y_{11,t}$, and in Ryde by $y_{12,t}$, 
also denote sexual offence counts in  Waverley by $y_{21,t}$, and in Ryde by $y_{22,t}$ ($t \in\{1,...,277\}$). 
As is recorded by the NSW Bureau of Crime Statistics and Research, the sexual offences are the sum of two subcategories:
(a) sexual assault; and (b) indecent assault, act of indecency and other sexual offences.
Thus at each time index $t$, the observation forms a $2\times2$ dimensional matrix-valued observation, denoted by $\bm{Y}_t$. 
Among these series, the first $t\in\{1,...,241\}$ counts form the training set used to fit the model, and the last three years $t\in\{242,...,277\}$ form the testing set, which serve as the real values for $h$-step ahead out-of-sample predictions.

\begin{figure}[!h]%
\centering%
\includegraphics[width=4in]{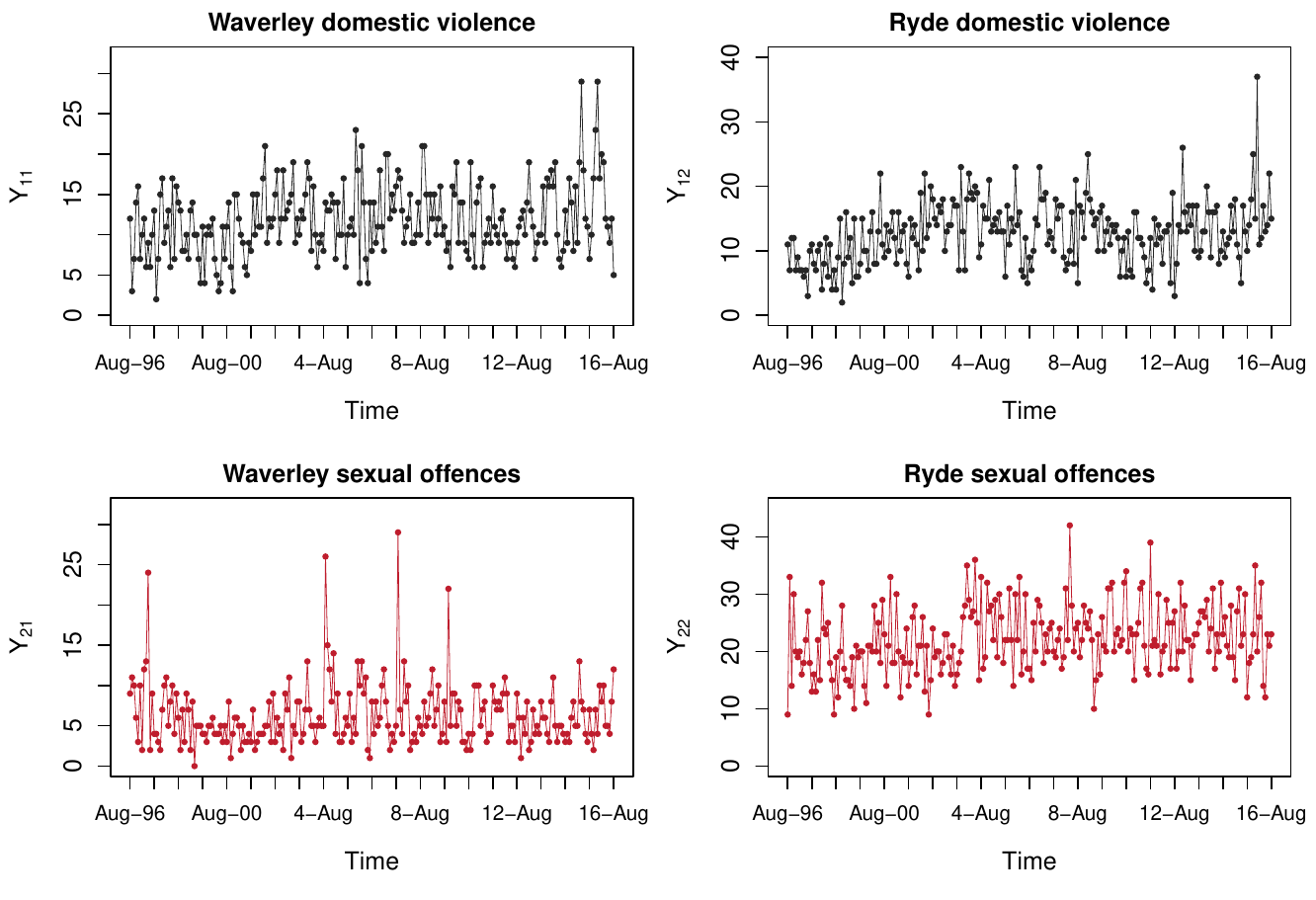}%
\caption{Time series of $\{\bm{Y}_t\}$  from August 1996 to August 2016.}
\label{trace1}
\end{figure}
\begin{figure}[!h]
\centering
\includegraphics[width=6in]{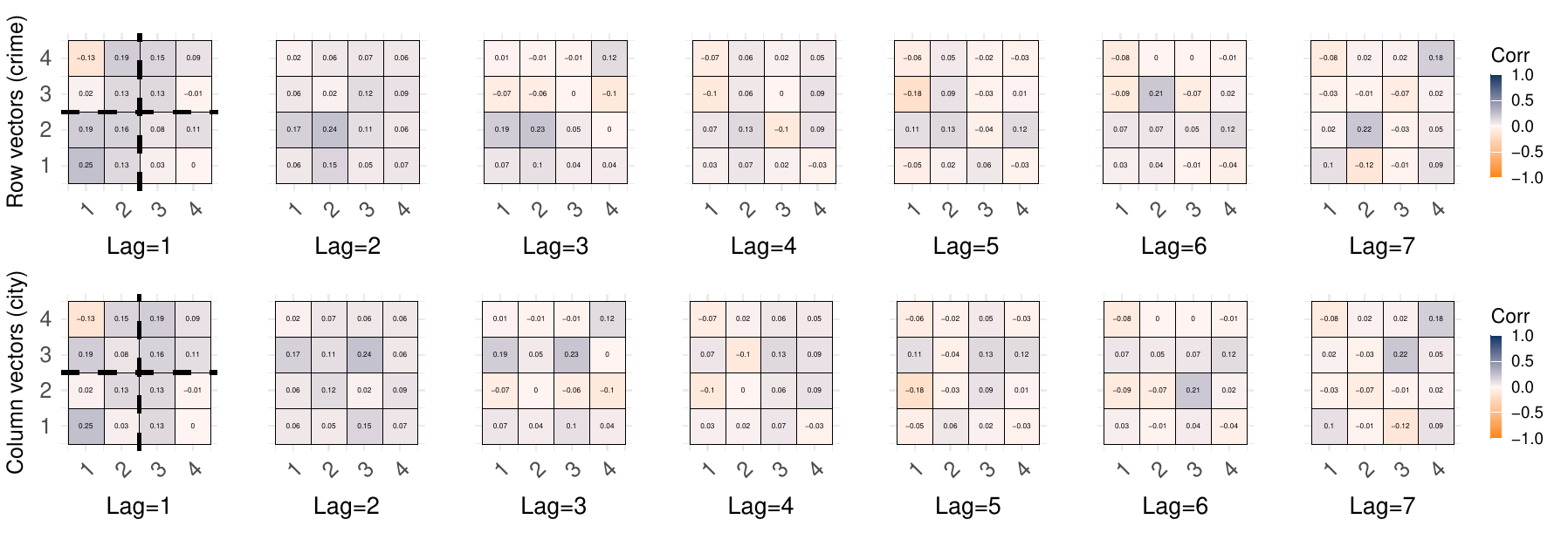}%
\caption{Cross-ACFs of $\{\bm{Y}_t\}$  from August 1996 to August 2016.}
\label{acf1}
\end{figure}

Figure \ref{trace1} shows the time series plots of $\{\bm{Y}_t\}$.  
As  seen in Figure \ref{trace1} that there is no clear trend in all series, indicating that all analyzed series are stationary. 
Moreover, we can calculate to see that the dispersion characteristics of each series are not the same, which 
brings challenges for model fitting. 
Figure \ref{acf1} shows the cross-ACF plots of $\{\bm{Y}_t\}$. 
Its top panel denotes the cross-ACFs of the same crime type across different cities, 
while 
the bottom panel shows the cross-ACFs of different criminal activities in a city. 
From the top panel of Figure \ref{acf1} we can see that all lags contain different  degrees of correlations,   
which suggests that the incidents of sexual offences (or domestic violence) in Waverley and Ryde are related. 
Moreover, 
the sexual offences in Waverly and Ryde are cross-correlated with domestic violence attacks in both cities. 
Similarly, the bottom panel shows us that Waverley's  domestic violence assults and sexual offences are related, and also interact with Ryde's two types of criminal incidents. 
Therefore, Figure \ref{acf1} shows that there are complex interdependencies between different cities, different types of crime, or across cities and types of crime. 
Therefore, it is suitable for matrix-variate modelling.

\begin{figure}[h]
\centering
\includegraphics[width=4in]{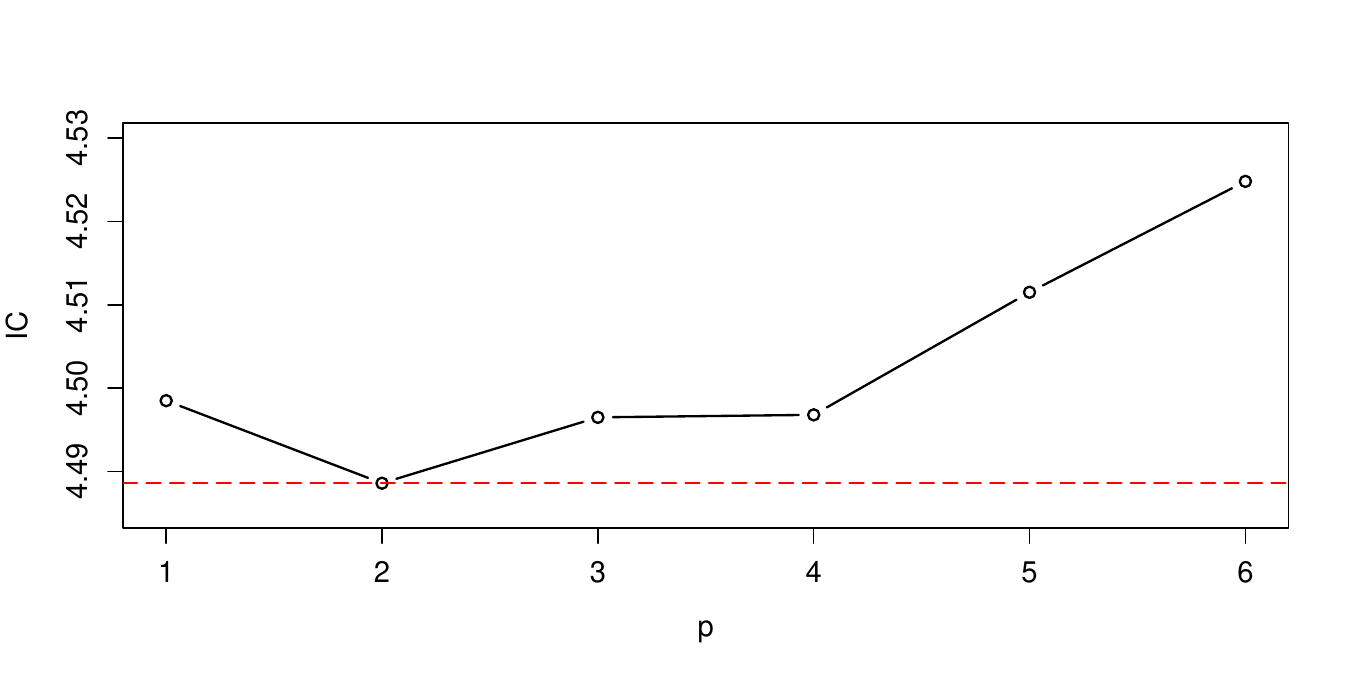}
\caption{The values of IC with different choice $p$. The horizontal red dashed line indicates the minimum $IC_1$ (\ref{IC1}) value.}
\label{c_p}
\end{figure}

Next, we use the MAT-INAR($p$) model to fit this data set, the order $p$ is chosen from 1 to 6. 
We use criterion (\ref{IC1}) to select the best order. 
To this end, we calculated $\hat{p}$ according to (\ref{IC1}) with $\overline{p}=6$, and draw 
the values of $IC_1$ in Figure \ref{c_p}. 
It can be seen from  Figure \ref{c_p} that a MAT-INAR(2) model is more suitable to fit the data. 
For comparative reasons,
we further choose the following models to fit this data set. 
The optimal orders of these vector models are calculated via criterion (E1) defined in Section E of Supplement Material:
\begin{itemize}
\item the stacked vector model, that is, MGINAR(4)  (\cite{Latour1997}) model;
\item two independent bivariate MINAR models (\cite{Pedeli2013a}) to fit domestic violence assaults  of Waverley and Ryde (MINAR(3)), and sexual offences of Waverley and Ryde (MINAR(1)), respectively (denoted by MINAR$^2_{\text{row}}$(3,1) model);
\item two independent bivariate MINAR models to fit sexual offences and domestic violence assaults  of Waverley (MINAR(1)) and of Ryde (MINAR(3)), respectively (denoted by MINAR$^2_{\text{col}}$(1,3) model);
\item the continuous matrix-variate autoregressive  (MAR(4))  model (\cite{Chen2021}); 
\item the MAT-INAR(2) model defined by (\ref{Mymodelp}) using the ICLS estimation method.
\end{itemize}
To select the best model among the competing models, we adopt  the
mean of residual sum of squares (MRSS) of each model.
That is, we fit the corresponding
models for training set and obtain the conditional expectation $E(\bm{Y}_{t}\vert \bm{Y}_{t-1})$ defined in (\ref{e1}).
The MRSS of each fitting model is given as follows
$
\text{MRSS}:={T}^{-1}\sum_{t=1}^T{\Big\Vert \bm{Y}_{t}-E(\bm{Y}_{t}\vert \bm{Y}_{t-1})\Big\Vert}_F.
$
We also obtain out-sample  forecast performances of  all  models for comparison purpose. 
Here we predict the occurrence of sexual offence and domestic violence assault of Waverley and Ryde in the next three years.
To illustrate the prediction effect, we consider the mean of the  out-of-sample prediction error (MSPE) as follows
$
\text{MSPE}={H}^{-1}\sum_{h=1}^H{\left\Vert \widehat{\bm Y}_{t}(h)-\bm{Y}_{t+h}\right\Vert}_F,
$
where  $H=36$, $\bm{Y}_{t+h}$ denotes the true observation with $t=241$, 
$\widehat{\bm{Y}}_{t}(h)$ is the correspondingly predictive value obtained via the $h$-step ahead conditional expectation. 
All the fitting and predicted results are summarized  in Table \ref{model_pk1}. 
We also give the number of parameters of the corresponding fitted model, expressed in $K$.

\begin{table}[hb]
\renewcommand\arraystretch{0.7}
\caption{Comparison between different models: MRSS, MSPE and $K$}
\label{model_pk1}
\centering
{\tabcolsep0.4in
\begin{tabular}{cccc}
  \hline
Model               & MRSS   & MSPE&$K$     \\
  \hline
 MAT-INAR(2)   &8.576   &12.589  &20  \\
 MGINAR(4)     &8.297   &12.669  &68 \\
MINAR$^2_{\text{row}}$(3,1)    &19.863    &36.342  &20     \\
MINAR$^2_{\text{col}}$(1,3)&11.294   &16.435  &20  \\
 MAR(4)   &8.705   &14.798  &36   \\
\hline
\end{tabular}
}
\end{table}

We analyze the results in Table \ref{model_pk1} in two ways.
On the surface, the MGINAR(4) model has the smallest
MRSS  among all models, while the MAT-INAR(2) model is only in a second place.  
However, the MGINAR(4) model has total 68 parameters, which 
are more than triple as many as the MAT-INAR(2) model with only 20 parameters. 
Meanwhile, the out-sample  forecasting performance of MGINAR(4) model is inferior to the MAT-INAR(2) model. 
Therefore, we conclude that 
a larger number of parameters result in its smaller MRSS, indicating the overfitting fact of the MGINAR(4) model. 
It can be seen that MAR(4) is not as good as MAT-INAR(2) and MGINAR(4) models in terms of fitting and prediction, which indicates that the continuous model has certain limitations for integer-valued data. 
Furthermore, 
it can be seen from the  results that the fitting and predicted effects of MINAR$^2_{\text{row}}$(3,1)   and MINAR$^2_{\text{col}}$(1,3)  models are poor,
which also indicates the necessity of establishing MAT-INAR model for the MITS.

\begin{figure}[!h]
\centering%
\includegraphics[width=4.5in]{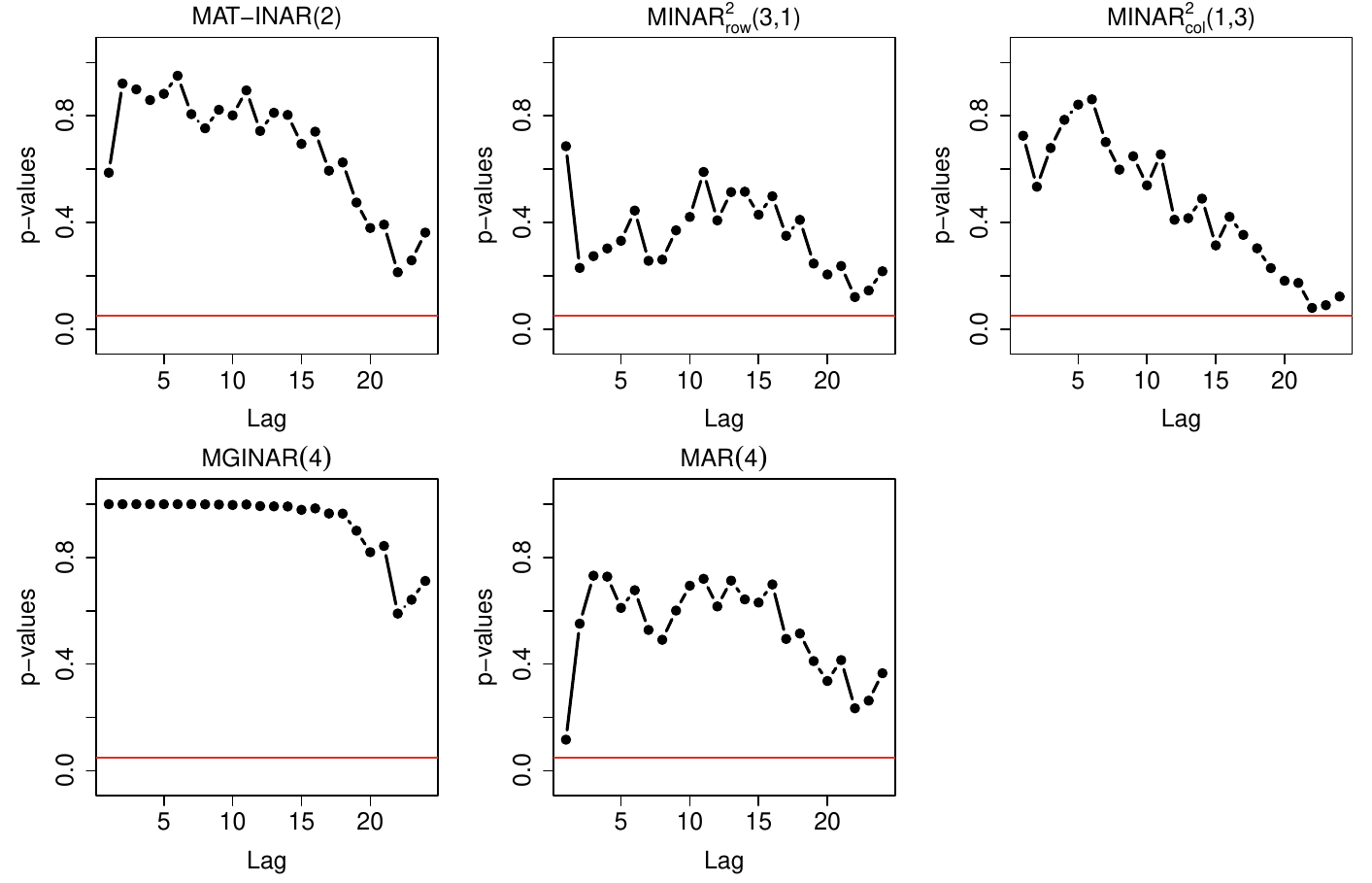}%
\caption{$p$-values of  portmanteau test statistics for residuals being serially uncorrelated. The red line indicates the significance level of 0.05.
}
\label{check}
\end{figure}

More intuitively, to assess the adequacy of  all different  models for the given dataset, it is natural to run some diagnostics
based on the residuals.
Since  MAT-INAR model can be converted to a vector model, standard diagnostics for vector models can be applied  (\cite{Chen2021,Tsay2024}).
Thus, we plot the $p$-values of portmanteau test with delay orders from 1
to 24 in Figure \ref{check}, which be applied to test for serial correlations among
the residual matrices.
We can see in Figure \ref{check} that the $p$-values of the portmanteau test  are all larger than 0.05 for all fitted models.
This also shows the effectiveness of the criteria (\ref{IC1}) and (E1), which can select the best model to fit the data adequately.

Now, we summarize the estimated coefficient matrix results of the MAT-INAR(2) model using the  ICLS method in Tables \ref{model_A}--\ref{model_V}, as well as the corresponding SE (in the parentheses) of $\widehat{\bm{A}}_{l,ICLS}$, $\widehat{\bm{B}}_{l,ICLS}~(l\in\{1,2\})$ and $\widehat{\bm{\Lambda}}_{ICLS}$. 
The significance of each parameter at $1\%$ (denoted $^{***}$), $5\%$ (denoted $^{**}$) and $10\%$ (denoted $^{*}$) level are 
also indicated, respectively.

The left and right autoregressive coefficient matrices results give some attractive conclusions. 
Table \ref{model_A} shows the $\widehat{\bm{A}}_{l,ils}~(l\in\{1,2\})$ results. 
For example, the first columns in $\widehat{\bm{A}}_{1,ils}$ and $\widehat{\bm{A}}_{2,ils}$
show the influence on the current sexual offence from the past  domestic violence assault.
The $\widehat{\bm{A}}_{1,ils}$ results indicate that the significant  dependence between domestic violence assault and  sexual offence. 
This is intuitively true, according to the surveys from the Australian Bureau of Statistics. 
Following the publication of an analysis on the SBS website in February 2020 entitled ``Australian Bureau of Statistics: Domestic violence is the number one cause of death for young and middle-aged women in Australia'', 
on average, one or more women die every week in Australia as a result of domestic or sexual violence, most of which are committed by 
close family members. 
Domestic violence often involves a prolonged form of power control and emotional abuse, 
making sexual assault more likely. 
As a result, domestic violence often sets the stage for sexual assault. 
While sexual assault is usually a one-time or short-term event, it does not necessarily change the dynamics of domestic violence or trigger the occurrence of domestic violence. 
Sexual offence is more likely to be a part of domestic violence, instead of adversely affecting the occurrence of domestic violence.

\begin{table}[h]
\renewcommand\arraystretch{0.8}
\caption{ Estimated left coefficient matrices ${\widehat{\bm{A}}}_{1,ils}$ and ${\widehat{\bm{A}}}_{2,ils}$ of MAT-INAR(2) model.}
\label{model_A}
\centering
{\tabcolsep0.02in
\begin{tabular}{cccccc}
  \hline
 \multirow{2}*{Estimation} &\multicolumn{2}{c}{${\widehat{\bm{A}}}_{1,ils}$}& & \multicolumn{2}{c}{${\widehat{\bm{A}}}_{2,ils}$}\\
 \cline{2-3}\cline{5-6}
              &Domestic violence &Sexual offence&&Domestic violence &Sexual offence\\
  \hline
 Domestic violence   &0.381 (0.055)$^{***}$ &0.047 (0.095)& &0.478 (0.048)$^{***}$ & 0.054 (0.114)\\
Sexual offence  &0.173 (0.080)$^{**}$  &0.254 (0.082)$^{***}$  & &0.121 (0.147)  & 0.117 (0.123)\\
\hline
\end{tabular}
}
\end{table}

\begin{table}[ht]
\renewcommand\arraystretch{0.8}
\caption{Estimated right coefficient matrices ${\widehat{\bm{B}}}_{1,ils}$ and ${\widehat{\bm{B}}}_{2,ils}$ of MAT-INAR(2) model.}
\label{model_B}
\centering
{\tabcolsep0.08in
\begin{tabular}{cccccc}
 \hline
 \multirow{2}*{Estimation} &\multicolumn{2}{c}{ ${\widehat{\bm{B}}}_{1,ils}$}& & \multicolumn{2}{c}{${\widehat{\bm{B}}}_{2,ils}$}\\
 \cline{2-3}\cline{5-6}
                 &Waverley &Ryde&   &Waverley &Ryde  \\
  \hline
 Waverley        &0.503 (0.144)$^{***}$ &0.084 (0.112)& &0.073 (0.125) & 0.238 (0.119)$^{**}$\\
  Ryde              &0.343 (0.196) $^{*}$ &0.393 (0.139)$^{***}$ &   &0.258 (0.153)$^{*}$ & 0.340 (0.130)$^{***}$\\
\hline
\end{tabular}
}
\end{table}

\begin{table}[hb]
\renewcommand\arraystretch{0.8}
\caption{ Estimated matrix $\widehat{\bm \Lambda}_{ils}$ of MAT-INAR(2) model.}
\label{model_V}
\centering
{\tabcolsep0.1in
\begin{tabular}{ccc}
 \hline
 \multirow{2}*{Estimation}  & \multicolumn{2}{c}{$\widehat{\bm \Lambda}_{ils}$}      \\
 \cline{2-3}
              &Waverley & Ryde\\
  \hline
 Domestic violence                &6.789 (1.224)$^{***}$    & 4.768 (1.936)$^{**}$  \\
Sexual offence                    &2.653 (1.289)$^{**}$            & 15.894 (1.957)$^{***}$ \\
\hline
\end{tabular}
}
\end{table}

Table \ref{model_B} shows the estimated $\bm{B}_l~(l\in\{1,2\})$ and  their effect should be considered in the view of $\bm{B}_l\circ_L\bm{Y}_{t-l}^{\top}$. 
It can be seen that the occurrence of crime offences between the two cities has a very significant cross-dependency, 
and the influence of Waverley on the occurrence of crime offences in Ryde is higher than that of Ryde on Waverley. 
We give some interpretations for $\widehat{\bm{B}}_{l,ils}$ results from a practical point of view. 
Given the geographical location of the two cities in Figure G1, these conclusions are meaningful.
As seen in Figure G1, we  can learn that Ryde and Waverley are on the south and north sides of the Parramatta River.
Therefore, geographical proximity is the main reason that Ryde is strongly correlated with Waverley.
Waverley is located in the east of Sydney, close to the coastline, and includes the famous Bondi Beach, is an area with a high concentration of tourist and commercial activities. 
Ryde, on the other hand, is located in the northwest of Sydney, far from the city centre. 
Violent crime in Waverley can quickly affect the surrounding area through the transport network (e.g. public transport, major roads, etc.), especially Ryde, which is easily accessible. 
Therefore, geographically, crime in Waverley is more likely to spread to Ryde, whereas crime in Ryde is  likely to less affect Waverley, as crime in Ryde may not directly affect more densely populated areas as in Waverley.

\begin{figure}[h]
\centering
\includegraphics[width=4in]{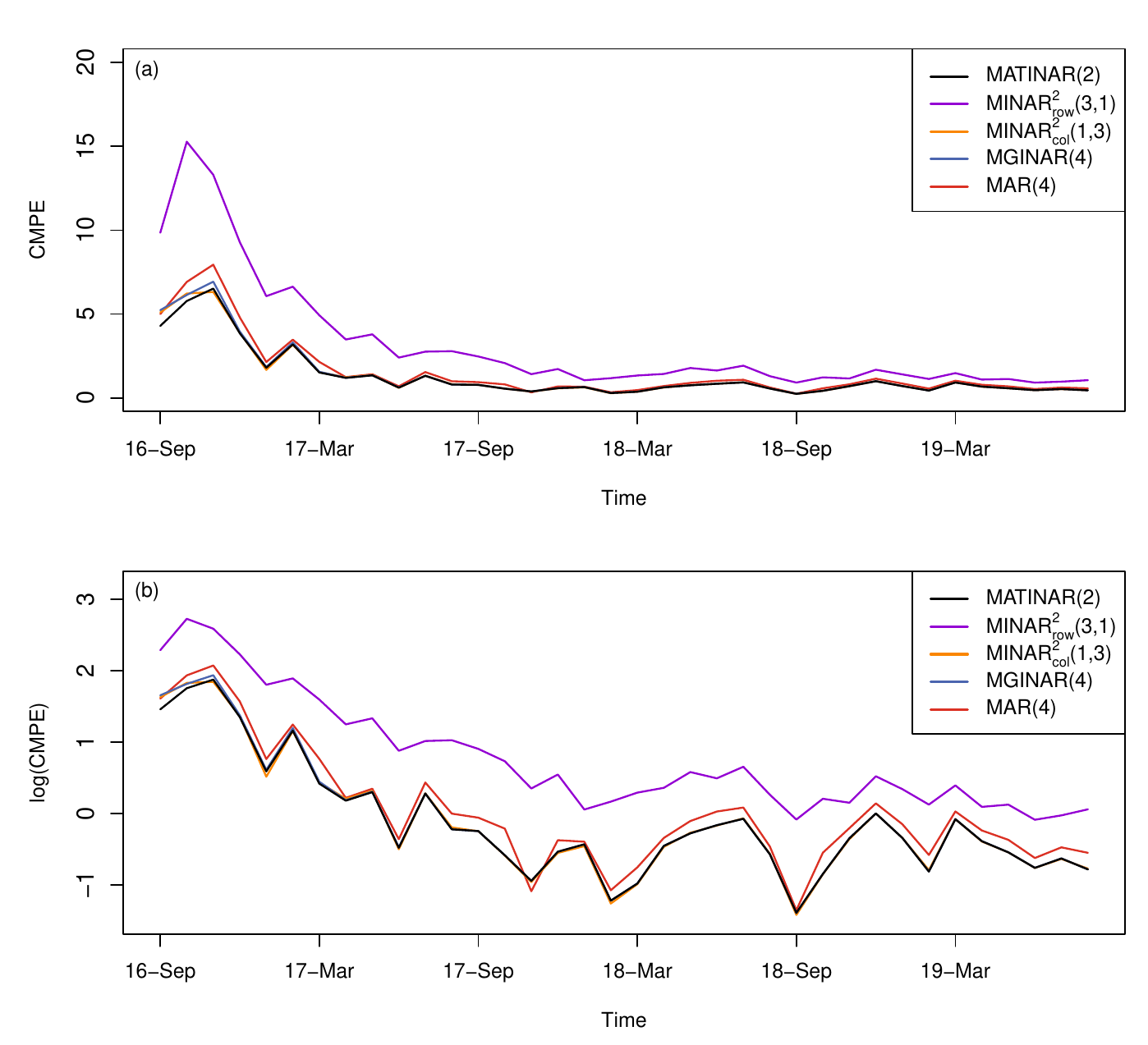}
\caption{CMPE curves  of all competitive models.
}
\label{CMPE_curve}
\end{figure}

We obtain the estimated $\widehat{\bm V}$ in Table \ref{model_V}.
Especially, we can see in Figure \ref{trace1},  the occurrence of  the sexual offence count in Ryde is serious. 
This is because Ryde has a larger city size and a larger population than Waverley, 
thus contributing to the occurance of sexual offence. 
Beside, cases of domestic violence and sexual offence may occur more frequently in economically disadvantaged Ryde, 
as poverty, unemployment and stress are often catalysts for domestic violence and sexual offence.

Finally, we focus  on the predictive performance of the fitted  models. 
For the sake of comparison,  we consider the cumulative loss function given by the mean-square predictive error (CMPE)   as follows 
$
{\text {CMPE}}_{S}={S}^{-1}\sum_{h=1}^{S}{\left\Vert \widehat{\bm Y}_{t}(h)-\bm{Y}_{t+h}\right\Vert}_F$ $(S\in\{1,...,36\}),
$
where $\bm{Y}_{t+h}$ denotes the true observation with $t=241$, 
$\widehat{\bm{Y}}_{t}(h)$ is the correspondingly predictive value obtained via the $h$-step ahead conditional expectation. 
The CMPE and the logarithmic CMPE curves are presented in Figure \ref{CMPE_curve}. 
We can see that the predictive power of MINAR$^2_{\text{row}}$(3,1) and MAR(4) are significantly worse.  
We can see that the CMPE values of MAT-INAR(2) model are mostly smaller than those of other models.
The MAT-INAR(2) model is more ideal for the prediction of MITS.
However, with the increase of predicted step size $S$, the CMPE curves of all models show a downward trend.
The fitting and predicted effects of MAT-INAR(2) model for each sequence of matrix-variate data are given in Figure G2 of Supplement Material, 
the right side of the vertical black dotted line is the $h$-step ahead out-of-sample prediction effect of all comparative models with $h\in\{1,...,6\}$.

\section{Conclusions}\label{sec:8}
\indent
In this study, we define two new matricial thinning operators and give some related properties. On this basis, 
the MAT-INAR($p$) model is constructed and the corresponding probabilistic and statistical properties are given. 
We give two estimation methods and the corresponding asymptotic theory. 
A new order selection criterion  is introduced  to address the order-determination problem.  We find that the MAT-INAR($p$) model achieves a substantial dimensional reduction by utilizing the matrix structure compared with the traditional vector model. 
The MAT-INAR($p$) model can maintain the original structure of matrix variate, reduce the information loss, and give a more realistic interpretation from the column and row variables. Besides,  for integer-valued data, the continuous model is no longer suitable, its fitting and prediction results are relatively poor. 
This also shows the necessity of constructing a  MAT-INAR model to provide a theoretical reference for MITS.

 There are a number of directions which are worth further investigations. On the one hand, the matrix-variate discrete distribution of more complex error terms can be considered to better capture the complexity and dispersion characteristics of the matrix-valued time series data. On the other hand, inspired by \cite{Zhang2024}, we may consider a multi-order additive matrix-variate integer-valued autoregressive model, which captures the row and column properties of $\bm{X}_{t-l}~(l\in\{1,...,p\})$ respectively.

\bibliographystyle{agsm}
\bibliography{mybibfile}
\newpage

{\bf Remark:}
This article is an improvement of the fifth chapter of the first author's master's thesis. The master's thesis was completed in June 2024. Currently, this article has been submitted and is under review.
\begin{figure}[!h]
\centering
\includegraphics[width=5in]{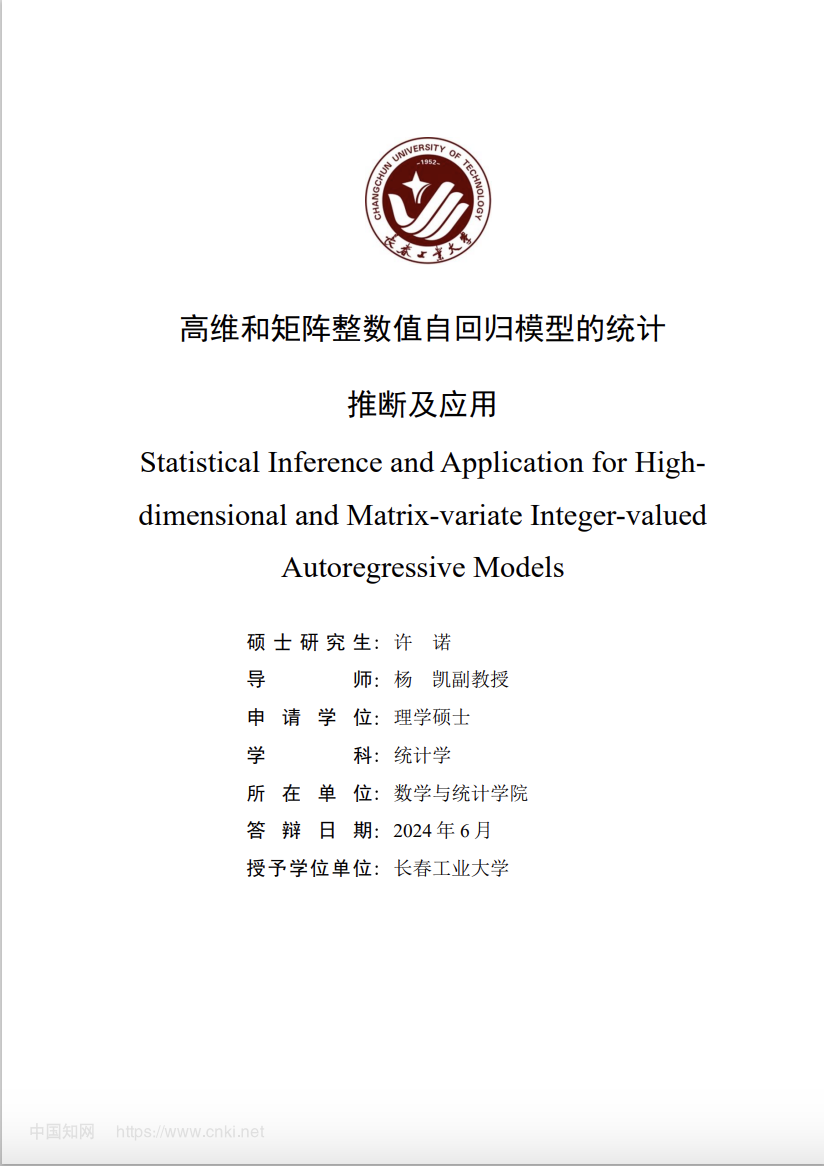}%
\end{figure}

\begin{figure}[!h]
\centering
\includegraphics[width=5in]{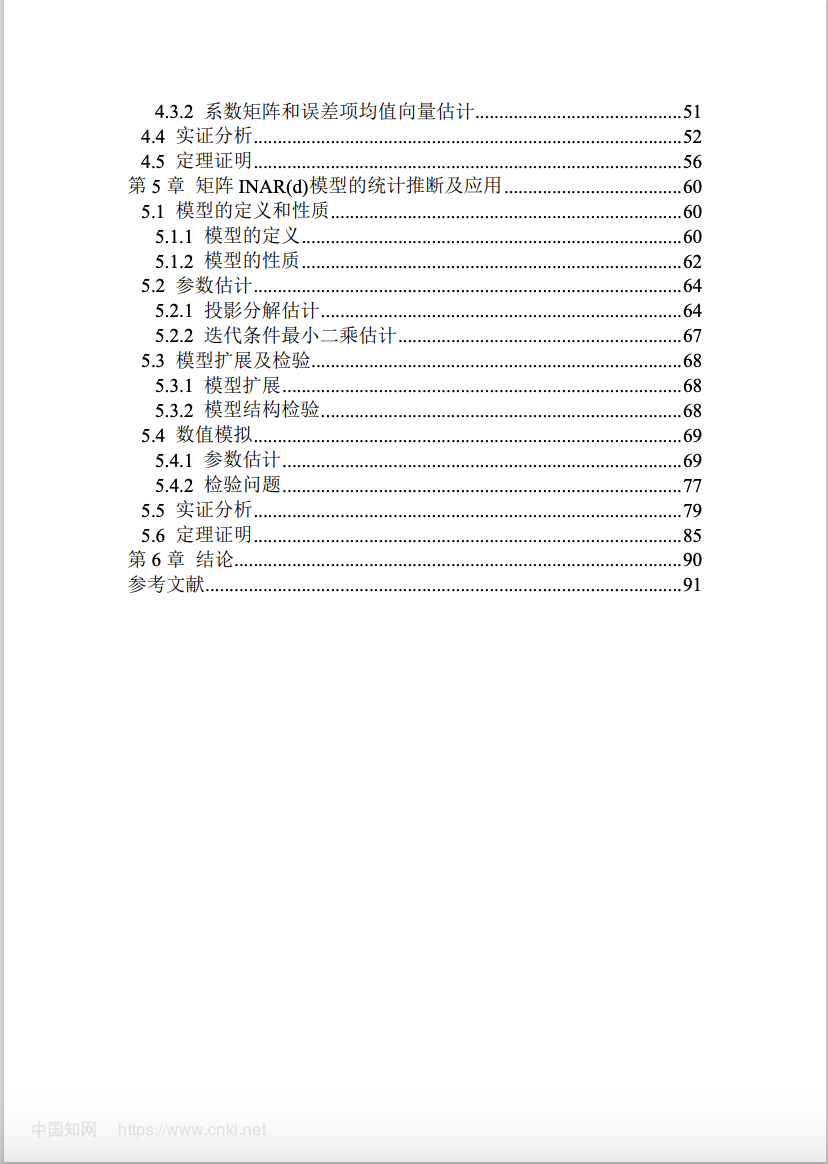}%
\end{figure}
\begin{figure}[!h]
\centering
\includegraphics[width=5.5in]{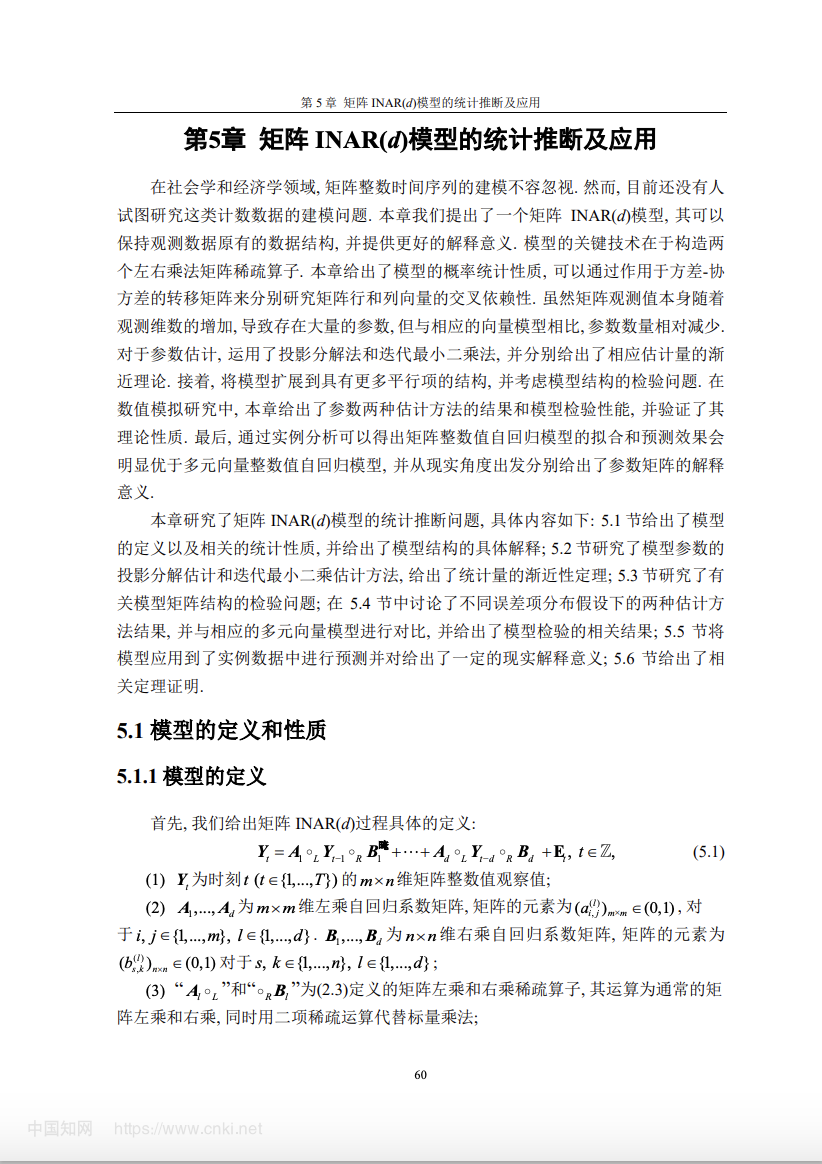}%
\end{figure}

\begin{figure}[!h]
\centering
\includegraphics[width=6in]{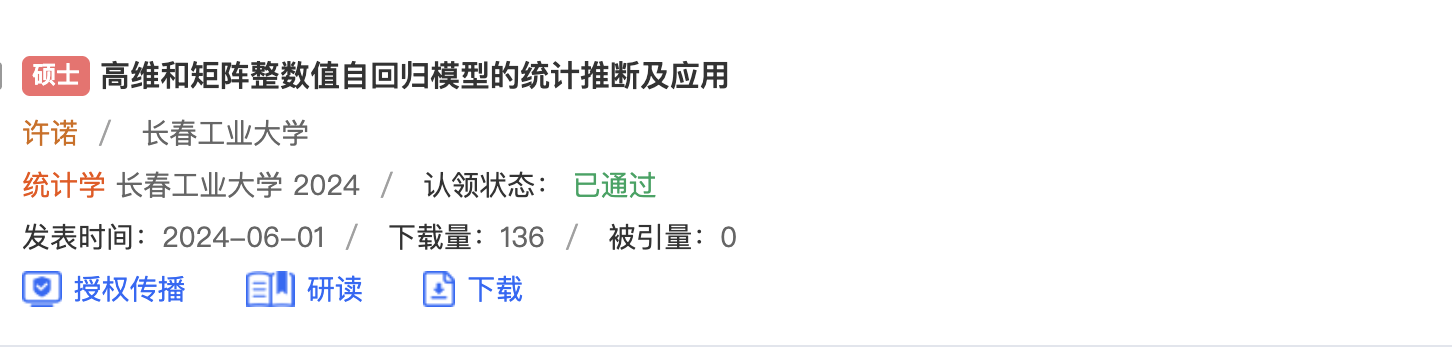}%
\end{figure}

\begin{figure}[!h]
\centering
\includegraphics[width=6in]{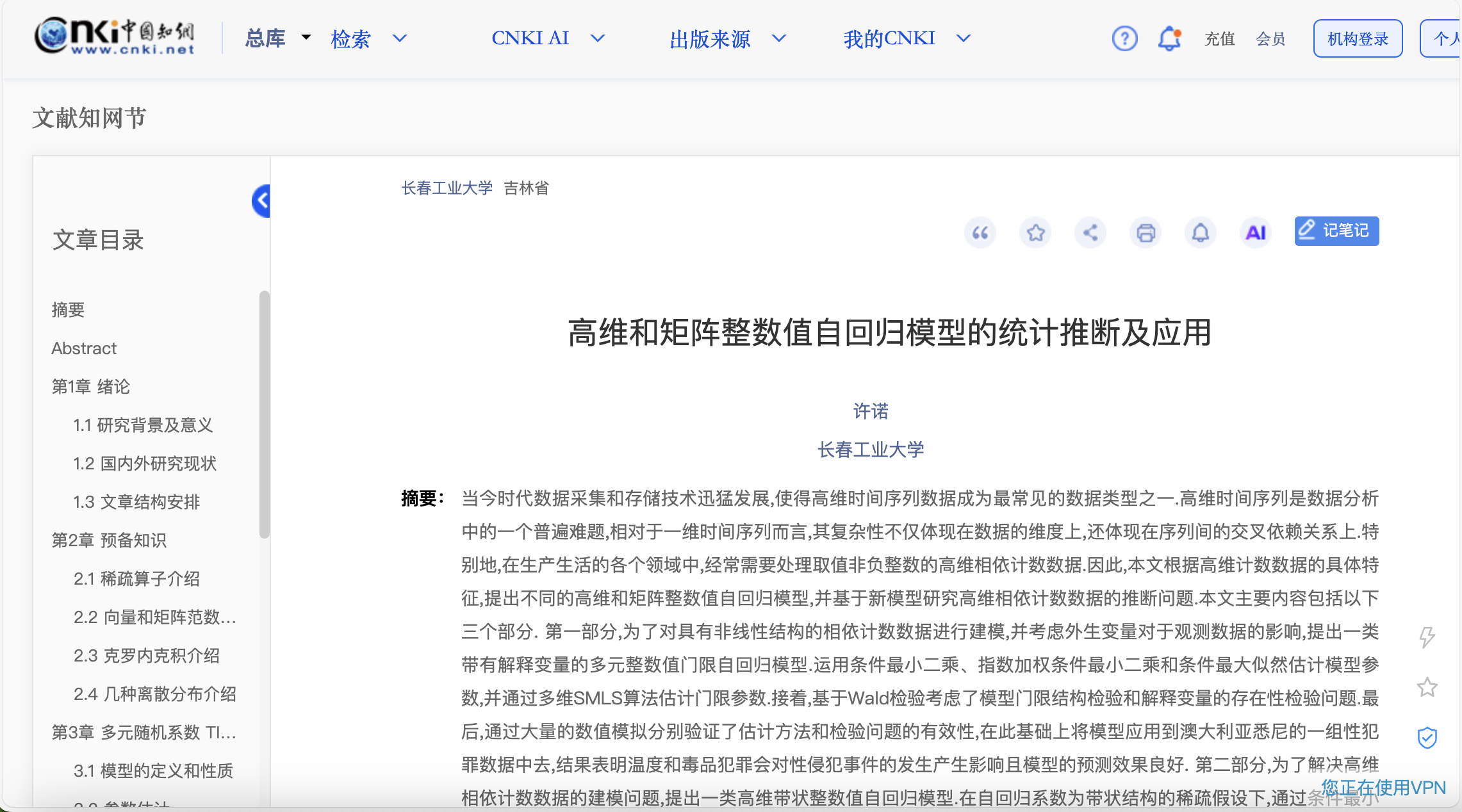}%
\end{figure}
\begin{figure}[!h]
\centering
\includegraphics[width=6in]{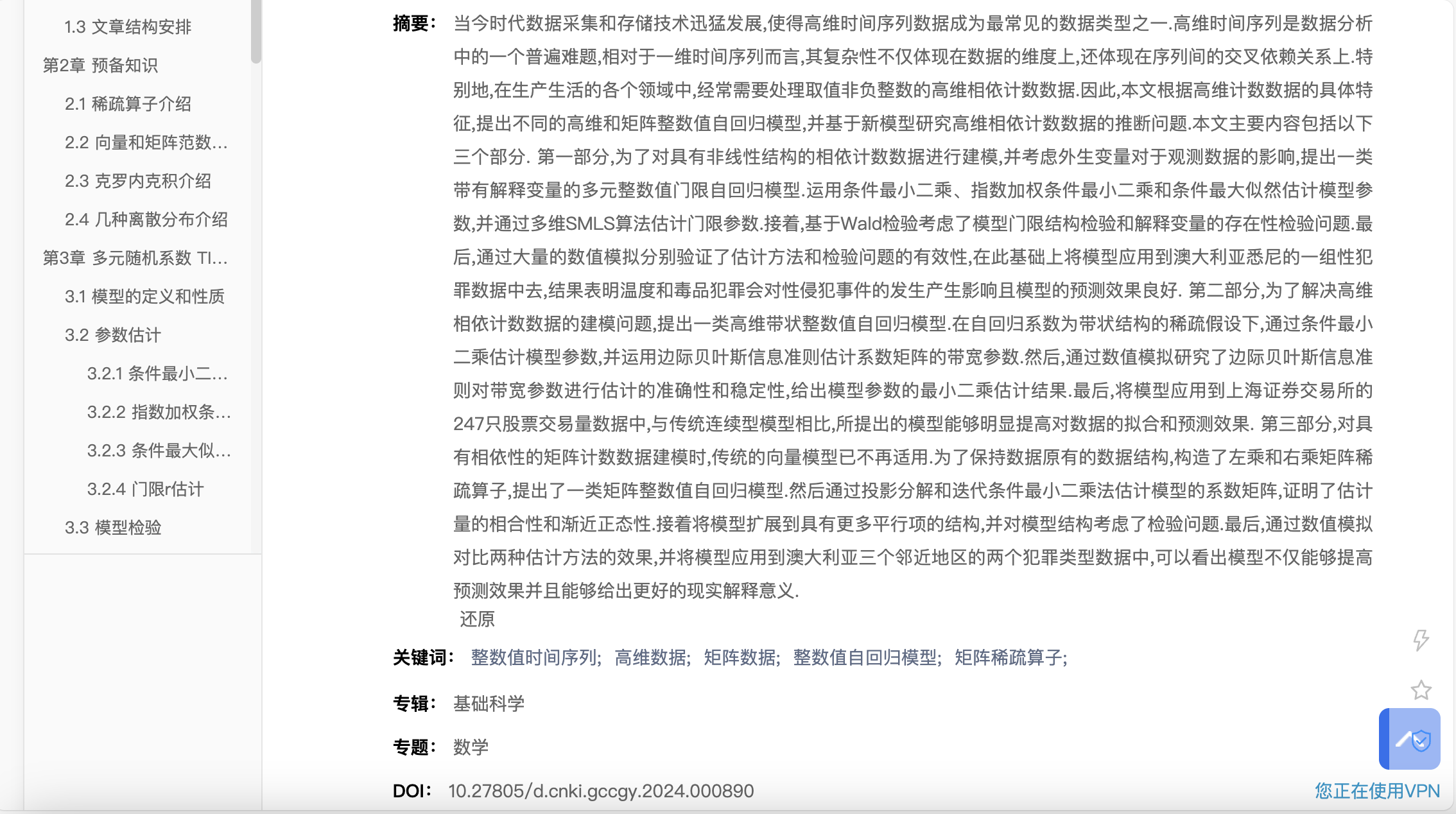}%
\end{figure}

\end{document}